\documentclass{amsart}
\input{diagrams}

\newtheorem{proposition}{Proposition}[section]
\newtheorem{lemma}[proposition]{Lemma}
\newtheorem{corollary}[proposition]{Corollary}
\newtheorem{theorem}[proposition]{Theorem}

\theoremstyle{definition}
\newtheorem{definition}[proposition]{Definition}

\theoremstyle{remark}

\def\yd{\mbox{$_H^H{\mathcal YD}$}}
\def\va{\varepsilon}
\def\v{\varphi}

\def\tr{\triangleright}

\def\ra{\rightarrow}
\def\a{\alpha}
\def\b{\beta}

\def\l{\lambda}
\def\r{\rho}
\def\cd{\cdot}

\def\ov{\overline}
\def\un{\underline}

\def\End{{\rm End}}

\newcommand{\smi}{\mbox{$S^{-1}$}}
\def\rawo\lonra{\longrightarrow}

\def\ot{\otimes}

\newcommand{\thlabel}[1]{\label{th:#1}}

\newcommand{\selabel}[1]{\label{se:#1}}
\newcommand{\seref}[1]{Section~\ref{se:#1}}
\newcommand{\lelabel}[1]{\label{le:#1}}
\newcommand{\leref}[1]{Lemma~\ref{le:#1}}
\newcommand{\prlabel}[1]{\label{pr:#1}}
\newcommand{\prref}[1]{Proposition~\ref{pr:#1}}
\newcommand{\colabel}[1]{\label{co:#1}}

\newcommand{\eqlabel}[1]{\label{eq:#1}}
\newcommand{\equref}[1]{(\ref{eq:#1})}

\begin{document}
\title[Yetter-Drinfeld modules over quasi-Hopf algebras]{More
properties of Yetter-Drinfeld modules over quasi-Hopf algebras}
\author{D. Bulacu}
\address{Faculty of Mathematics, University of Bucharest,
RO-70109 Bucharest 1, Romania}
\email{dbulacu@al.math.unibuc.ro}

\author{S. Caenepeel}
\address{Faculty of Applied Sciences,
Vrije Universiteit Brussel, VUB, B-1050 Brussels, Belgium}
\email{scaenepe@vub.ac.be}
\urladdr{http://homepages.vub.ac.be/\~{}scaenepe/}
\author{F. Panaite}
\address{Institute of Mathematics of the Romanian
Academy,
PO-Box 1-764, RO-70700 Bucharest, Romania}
\email{florin.panaite@imar.ro}
\thanks{This research supported by the bilateral project
``Hopf Algebras in Algebra, Topology, Geometry and Physics" of the Flemish and
Romanian governments. The third author was also partially supported by
the programmes SCOPES and EURROMMAT. The first and third author wish
to the Vrije Universiteit Brussel for its warm hospitality during
their visit there.}
\begin{abstract}
We generalize various properties of Yetter-Drinfeld modules over
Hopf algebras to quasi-Hopf algebras. The dual of a finite dimensional
Yetter-Drinfeld module is again a Yetter-Drinfeld module. The algebra $H_0$
in the category of Yetter-Drinfeld modules that can be obtained by
modifying the multiplication in a proper way is quantum commutative.
We give a Structure Theorem for Hopf modules in the category of Yetter-Drinfeld
modules, and deduce the existence and uniqueness of integrals from it.
\end{abstract}
\maketitle

\section{Introduction}\selabel{1}
The motivation for studying Yetter-Drinfeld modules over quasi-Hopf algebras
is the same as for Hopf algebras. It is well known that for
any finite dimensional
Hopf algebra $H$ the category of Yetter-Drinfeld modules ${}_H{\mathcal YD}^H$
is isomorphic to the category of modules
over the quantum double $D(H)$. From a categorical point of view, the quantum
double $D(H)$ arises by considering the center ${\mathcal Z}({}_H{\mathcal M})$ of
the monoidal category ${}_H{\mathcal M}$ of left $H$-modules. More precisely, one
has ${\mathcal Z}({}_H{\mathcal M})\simeq {}_{D(H)}{\mathcal M}$ if $H$ is finite dimensional. Actually, the category of Yetter-Drinfeld modules appears as
an intermediate step in the proof of this isomorphism: one first proves that
${\mathcal Z}({}_H{\mathcal M})\simeq {}_H{\mathcal YD}^H$,
and then ${}_H{\mathcal YD}^H\simeq {}_{D(H)}{\mathcal M}$, where the finite
dimensionality is not needed in the proof of the first isomorphism,
see \cite{k} for full detail.\\
Quasi-bialgebras and quasi-Hopf algebras were introduced by Drinfeld \cite{d1};
a categorical interpretation is the following: a quasi-bialgebra $H$
is an algebra with the additional structure that is needed to make the
category of left $H$-modules, with the tensor product over $k$ as tensor
product and $k$ as unit object into a monoidal category. The difference
with a usual bialgebra is that we do not require that the associativity
isomorphism coincides with the associativity in the category of vector spaces.
A quasi-Hopf algebra is a quasi-bialgebra with  additional structure
making the category of finite dimensional $H$-modules into a
monoidal category with duality.\\
The center construction ${\mathcal Z}({\mathcal C})$ can be applied to any monoidal category ${\mathcal C}$. Majid \cite{m1} computed the center of the category
of left modules over a quasi-Hopf algebra $H$, and introduced the
category of Yetter-Drinfeld modules over $H$. Hausser and Nill \cite{hn1}, \cite{hn2} constructed the quantum double $D(H)$ of a
finite dimensional quasi-Hopf algebra $H$, and proved that ${}_H{\mathcal YD}^H\simeq
{}_{D(H)}{\mathcal M}$. Recently, Schauenburg \cite{sch} gave the equivalence
between the category of Yetter-Drinfeld modules $\yd $ and the category
${}_H^H{\mathcal M}_H^H$ of Hopf bimodules.
In \cite{bn}, the relation between Yetter-Drinfeld modules and
Radford's biproduct is studied.
In \cite{bcp}, the rigidity of the category of Yetter-Drinfeld modules is
investigated, as well as the relations between left, left-right, right-left
and right Yetter-Drinfeld modules.\\
In this paper, which can be seen as a sequel to \cite{bcp}, we continue our
investigations of properties of Yetter-Drinfeld modules. In \seref{3},
we show that the linear dual of a finite dimensional right-left
Yetter-Drinfeld module is a left-right Yetter-Drinfeld module.\\
It was shown in
\cite{bpv}, \cite{bn} that the multiplication on $H$ can be modified in such
a way that we obtain an algebra in the category of left Yetter-Drinfeld modules.
The main result of \seref{4} is that $H_0$ is quantum commutative.\\
In \seref{5}, we will generalize Doi's results \cite{doi} about Hopf modules in
the category of Yetter-Drinfeld modules to our situation: we give
a Structure Theorem for Hopf modules in the category of Yetter-Drinfeld
modules over a quasi-Hopf algebras, and we use this result to obtain
the existence and uniqueness of integrals for a finite dimensional
braided Hopf algebra in $\yd $. We apply this to the braided Hopf algebra
considered in \seref{4}, in the case where $H$ is finite dimensional
and quasitriangular.

\section{Preliminary results}\selabel{2}
\subsection{Quasi-Hopf algebras}
We work over a commutative field $k$. All algebras, linear spaces
etc. will be over $k$; unadorned $\ot $ means $\ot_k$. Following
Drinfeld \cite{d1}, a quasi-bialgebra is a fourtuple $(H, \Delta ,
\va , \Phi )$, where $H$ is an associative algebra with unit,
$\Phi$ is an invertible element in $H\ot H\ot H$, and $\Delta :\
H\ra H\ot H$ and $\va :\ H\ra k$ are algebra homomorphisms
satisfying the identities
\begin{eqnarray}
&&(id \ot \Delta )(\Delta (h))=
\Phi (\Delta \ot id)(\Delta (h))\Phi ^{-1},\label{q1}\\
&&(id \ot \va )(\Delta (h))=h\ot 1, 
\mbox{${\;\;\;}$}
(\va \ot id)(\Delta (h))=1\ot h,\label{q2}
\end{eqnarray}
for all $h\in H$, and $\Phi$ has to be a normalized $3$-cocycle,
in the sense that
\begin{eqnarray}
&&(1\ot \Phi)(id\ot \Delta \ot id) (\Phi)(\Phi \ot 1)= (id\ot id
\ot \Delta )(\Phi ) (\Delta \ot id \ot id)(\Phi
),\label{q3}\\
&&(id \ot \va \ot id )(\Phi )=1\ot 1\ot 1.\label{q4}
\end{eqnarray}
The map $\Delta $ is called the coproduct or the
comultiplication, $\va $ the counit and $\Phi $ the reassociator.
As for Hopf algebras \cite{sw} we use the notation $\Delta (h)=\sum h_1\ot
h_2$. Since $\Delta$ is only quasi-coassociative we adopt the
further notation
$$
(\Delta \ot id)(\Delta (h))= \sum h_{(1, 1)}\ot h_{(1, 2)}\ot h_2,
\mbox{${\;\;\;}$} (id\ot \Delta )(\Delta (h))=
\sum h_1\ot h_{(2, 1)}\ot h_{(2,2)}, 
$$
for all $h\in H$. We will denote the tensor components of $\Phi$
by capital letters, and the ones of $\Phi^{-1}$ by small letters,
namely
\begin{eqnarray*}
&&\Phi=\sum X^1\ot X^2\ot X^3= \sum T^1\ot T^2\ot T^3= \sum V^1\ot
V^2\ot V^3=\cdots\\%
&&\Phi^{-1}=\sum x^1\ot x^2\ot x^3= \sum
t^1\ot t^2\ot t^3= \sum v^1\ot v^2\ot v^3=\cdots
\end{eqnarray*}
A quasi-bialgebra $H$ is called a quasi-Hopf algebra if there exists an
anti-automorphism $S$ of the algebra $H$ and $\a , \b \in
H$ such that:
\begin{eqnarray}
&&\sum S(h_1)\a h_2=\va (h)\a \mbox{${\;\;\;}$ and ${\;\;\;}$}
\sum h_1\b S(h_2)=\va (h)\b ,\label{q5}\\[1mm]%
&&\sum X^1\b S(X^2)\a X^3=1 %
\mbox{${\;\;\;}$ and${\;\;\;}$}%
\sum S(x^1)\a x^2\b S(x^3)=1,\label{q6}
\end{eqnarray}
for all $h\in H$. It is shown in  \cite{bc} that the condition that the antipode
is bijective follows automatically from the other axioms in the case where
$H$ is finite dimensional. Observe that
 the antipode of a quasi-Hopf algebra is determined uniquely up to a
transformation $\a \mapsto U\a $, $\b \mapsto \b
U^{-1}$, $S(h)\mapsto US(h)U^{-1}$, where $U\in H$ is invertible.
The axioms for a quasi-Hopf algebra imply that $\va (\a )\va (\b
)=1$, so, by rescaling $\a $ and $\b $, we may assume without loss
of generality that $\va (\a )=\va (\b )=1$ and $\va \circ S=\va $.
The identities (\ref{q2}-\ref{q4}) also imply
that
\begin{equation}\label{q7}
(\va \ot id\ot id)(\Phi )= (id \ot id\ot \va )(\Phi )=1\ot 1\ot 1.
\end{equation}
Together with a quasi-Hopf algebra 
$H=(H, \Delta , \va , \Phi , S, \a , \b )$ we also have $H^{\rm op}$, $H^{\rm cop}$
and $H^{\rm op, cop}$ as quasi-Hopf algebras, where ``op" means opposite
multiplication and ``cop" means opposite comultiplication. The 
reassociators of these three quasi-Hopf algebras are
$\Phi_{\rm op}=\Phi ^{-1}$,
$\Phi_{\rm cop}=(\Phi ^{-1})^{321}$, $\Phi_{\rm op, cop}=\Phi ^{321}$,
the antipodes are
$S_{\rm op}=S_{\rm cop}=(S_{\rm op,cop})^{-1}=S^{-1}$,
and the elements $\alpha,\beta$ are $\a_{\rm op}=\smi (\b )$,
$\b _{\rm op}=\smi (\a )$, $\a _{\rm cop}=\smi (\a )$, $\b _{\rm cop}=\smi (\b )$,
$\a _{\rm op,cop}=\b $ and $\b _{\rm op,cop}=\a $.\\
Recall next that the definition of a quasi-Hopf algebra is
``twist coinvariant", in the following sense. An invertible element
$F\in H\ot H$ is called a {\sl gauge transformation} or {\sl
twist} if $(\va \ot id)(F)=(id\ot \va)(F)=1$. If $H$ is a
quasi-Hopf algebra and $F=\sum F^1\ot F^2\in H\ot H$ is a gauge
transformation with inverse $F^{-1}=\sum G^1\ot G^2$, then we can
define a new quasi-Hopf algebra $H_F$ by keeping the
multiplication, unit, counit and antipode of $H$ and replacing the
comultiplication, antipode and the elements $\alpha$ and $\beta$
by
\begin{eqnarray}
&&\Delta _F(h)=F\Delta (h)F^{-1},\label{g1}\\[1mm]
&&\Phi_F=(1\ot F)(id \ot \Delta )(F) \Phi (\Delta \ot id)
(F^{-1})(F^{-1}\ot 1),\label{g2}\\
&&\a_F=\sum S(G^1)\a G^2,
\mbox{${\;\;\;}$}
\b_F=\sum F^1\b S(F^2).\label{g3}
\end{eqnarray}
It is well-known that the antipode of a Hopf
algebra is an anti-coalgebra morphism. The corresponding statement
for a quasi-Hopf algebra is the following: there exists a gauge
transformation $f\in H\ot H$ such that
\begin{equation} \label{ca}
f\Delta (S(h))f^{-1}= \sum (S\ot S)(\Delta ^{\rm cop}(h)),
\end{equation}
for all $h\in H$,
where $\Delta ^{\rm cop}(h)=\sum h_2\ot h_1$. The element
$f$ can be computed explicitly. First set
\begin{equation}
\sum A^1\ot A^2\ot A^3\ot A^4= (\Phi \ot 1) (\Delta \ot id\ot
id)(\Phi ^{-1}),
\end{equation}
\begin{equation} \sum B^1\ot B^2\ot B^3\ot B^4=
(\Delta \ot id\ot id)(\Phi )(\Phi ^{-1}\ot 1)
\end{equation}
and then define $\gamma, \delta\in H\ot H$ by
\begin{equation} \label{gd}
\gamma =\sum S(A^2)\a A^3\ot S(A^1)\a A^4~~{\rm and}~~ \delta
=\sum B^1\b S(B^4)\ot B^2\b S(B^3).
\end{equation}
Then $f$ and $f^{-1}$ are given by the formulas
\begin{eqnarray}
f&=&\sum (S\ot S)(\Delta ^{\rm op}(x^1)) \gamma \Delta (x^2\b
S(x^3)),\label{f}\\
f^{-1}&=&\sum \Delta (S(x^1)\a x^2) \delta (S\ot S)(\Delta
^{\rm op}(x^3)).\label{g}
\end{eqnarray}
Moreover, $f$ satisfies the following relations:
\begin{equation} \label{gdf}
f\Delta (\a )=\gamma , \mbox{${\;\;\;}$} 
\Delta (\b )f^{-1}=\delta .
\end{equation}
Furthermore the corresponding twisted reassociator (see
(\ref{g2})) is given by
\begin{equation} \label{pf}
\Phi _f=\sum (S\ot S\ot S)(X^3\ot X^2\ot X^1).
\end{equation}
In a Hopf algebra $H$, we obviously have the identity
$$\sum h_1\ot h_2S(h_3)=h\ot 1,~{\rm for~all~}h\in H.$$
We will need the generalization of this formula to
the quasi-Hopf algebra setting.
Following \cite{hn1,hn2}, we define
\begin{eqnarray} \label{qr}
p_R&=&\sum p^1_R\ot p^2_R=\sum x^1\ot x^2\b S(x^3),\\
q_R&=&\sum q^1_R\ot q^2_R=\sum X^1\ot S^{-1}(\a X^3)X^2,\label{qra}\\
\label{ql}
p_L&=&\sum p^1_L\ot p^2_L=\sum X^2\smi (X^1\b )\ot X^3,\\
q_L&=&\sum q^1_L\ot q^2_L=\sum S(x^1)\a x^2\ot x^3.\label{qla}
\end{eqnarray}
We then have, for all $h\in H$,
\begin{eqnarray} \label{qr1}
\sum \Delta (h_1)p_R[1\ot S(h_2)]&=&p_R(h\ot 1),\\
\sum [1\ot S^{-1}(h_2)]q_R\Delta (h_1)&=&(h\ot 1)q_R,\\
\label{ql1}
\sum \Delta (h_2)p_L[\smi (h_1)\ot 1]&=&p_L(1\ot h),\\
\sum [S(h_1)\ot 1]q_L\Delta (h_2)&=&(1\ot h)q_L,
\end{eqnarray}
and
\begin{eqnarray}
&&\hspace*{-3cm}(q_R\ot 1)(\Delta \ot id)(q_R)\Phi^{-1}=
\sum [1\ot S^{-1}(X^3)\ot S^{-1}(X^2)]\nonumber\\
&& [1\ot S^{-1}(f^2)\ot
S^{-1}(f^1)] (id \ot \Delta )(q_R\Delta(X^1)),\label{qr2}
\end{eqnarray}
where $f=\sum f^1\ot f^2$ is the twist defined in (\ref{f}).\\
A quasi-Hopf algebra $H$ is quasitriangular if
there exists an element $R\in H\ot H$ such that
\begin{eqnarray}
(\Delta \ot id)(R)&=&\sum \Phi _{312}R_{13}\Phi ^{-1}_{132}R_{23}\Phi ,\label{qt1}\\
(id \ot \Delta )(R)&=&\sum \Phi ^{-1}_{231}R_{13}\Phi _{213}R_{12}\Phi ^{-1},\label{qt2}\\
\Delta ^{\rm cop}(h)R&=&R\Delta (h),~{\rm for~all~}h\in H,\label{qt3}\\
(\va \ot id)(R)&=&(id\ot \va)(R)=1.\label{qt4}
\end{eqnarray}
Here we used the following notation:
if $\sigma $ is a permutation of $\{1, 2, 3\}$, then we write $\Phi _{\sigma (1)
\sigma (2)\sigma (3)}=\sum X^{\sigma ^{-1}(1)}\ot
X^{\sigma ^{-1}(2)}\ot X^{\sigma ^{-1}(3)}$;  $R_{ij}$ means $R$ acting non-trivially
on the $i$-{th} and $j$-{th} tensor factors of $H\ot H\ot H$.\\
It is shown in {\cite{bn3}} that $R$ is
invertible. Furthermore, the element
\begin{equation} \label{elmu}
u=\sum S(R^2p^2)\a R^1p^1,
\end{equation}
with $p_R=\sum p^1\ot p^2$ defined as in (\ref{qr}), is
invertible in $H$, and
\begin{equation} \label{inelmu}
u^{-1}=\sum X^1R^2p^2S(S(X^2R^1p^1)\a X^3),
\end{equation}
\begin{equation} \label{sqina}
\va (u)=1~~{\rm and}~~
S^2(h)=uhu^{-1},
\end{equation}
for all $h\in H$. Consequently the antipode $S$ is bijective, so,
as in the Hopf algebra case, the assumptions about invertibility
of $R$ and bijectivity of $S$ can be dropped. Moreover, the
$R$-matrix $R=\sum R^1\ot R^2$ satisfies the identity (see
\cite{ac}, \cite{hn2}, \cite{bn3}):
\begin{equation} \label{ext}
f_{21}Rf^{-1}=(S\ot S)(R)
\end{equation}
where $f=\sum f^1\ot f^2$ is the twist defined in (\ref{f}), and
$f_{21}=\sum f^2\ot f^1$.

\subsection{Monoidal categories}
A monoidal or tensor category is a sixtuple $({\mathcal C},\ot, \un{1},a,l,r)$,
where ${\mathcal C}$ is a category, $\ot$ is a functor
${\mathcal C}\times {\mathcal C}\ra {\mathcal C}$ (called the tensor product),
$\un{1}$ is an object of ${\mathcal C}$, and
$$a_{U, V, W}:\ (U\ot V)\ot W\ra U\ot (V\ot W)$$
$$l_V:\ V\cong V\ot \un{1}~~;~~r_V:\ V\cong \un{1}\ot V$$
are natural isomorphisms satisfying certain coherence conditions,
see for example \cite{k,McLane,m2}.
An object $V$ of a monoidal category ${\mathcal C}$
has a left dual if there exists
an object $V^*$ and morphisms ${\rm ev}_V:\ V^*\ot V\ra \un{1}$,
${\rm coev}_V:\ \un{1}\ra V\ot V^*$ in ${\mathcal C}$ such that
\begin{eqnarray}
&&\hspace*{-2cm}
l^{-1}_V\circ (id_V\ot ev_V)
\circ a_{V, V^*, V}\circ (coev_V\ot id_V)\circ r_V=id_V,\label{rig1}\\
&&\hspace*{-2cm}
r^{-1}_{V^*}\circ (ev_V\ot id_{V^*})
\circ a^{-1}_{V^*, V, V^*}\circ
(id_{V^*}\ot coev_V)\circ l_{V^*}=id_{V^*}\label{rig2}.
\end{eqnarray}
${\mathcal C}$ is called a rigid monoidal category if every object of
${\mathcal C}$ has a dual.\\
A braided monoidal category is a monoidal category
equipped with a commutativity natural isomorphism
$c_{U, V}:\ U\ot V\ra V\ot U$,
compatible with the unit and the associativity.\\
In a braided monoidal category, we can define algebras, coalgebras,
bialgebras and Hopf algebras. For example, a bialgebra $(B, \un{m}, \un{\eta }, \un {\Delta }, \un {\va })$ consists of $B\in {\mathcal C}$, a
multiplication $\un{m}:\ B\ot B\ra B$ which is associative up to the
natural isomorphism $a$, and a unit $\un {\eta }: \un {1}\ra B$ such that
$\un{m}\circ (\un{\eta }\ot id)=\un {m}\circ (id \ot
\un{\eta })=id$. The properties of the comultiplication $\un{\Delta}$
and the counit $\un {\va }$ are similar. In addition,
$\un{\Delta}:\ B\to B\ot B$ has to be an algebra morphism, where $B\ot B$
is an algebra with multiplication $\un{m}_{B\ot B}$, defined as the
composition
\begin{equation}\label{bi}
\begin{array}{ccc}
(B\ot B)\ot(B\ot B)&\rTo^{a}& B\ot(B\ot (B\ot B))\\
&\rTo^{id\ot a^{-1}}& B\ot((B\ot B)\ot B)\\
&\rTo^{id\ot c\ot id}&B\ot((B\ot B)\ot B)\\
&\rTo^{id\ot a}&B\ot(B\ot (B\ot B))\\
&\rTo^{a^{-1}}&(B\ot B)\ot(B\ot B)\\
&\rTo^{\un{m}\ot\un{m}}&B\ot B
\end{array}
\end{equation}
A Hopf algebra $B$ is a bialgebra with a
morphism $\un {S}: B\ra B$ in ${\mathcal C}$ (the antipode)
satisfying the usual axioms
$\un {m}\circ (\un {S}\ot id)\circ \un {\Delta }=\un {\eta }\circ \un {\va }=
\un {m}\circ (id \ot \un {S})\circ \un {\Delta }$.
It is known, see e.g. \cite{m3},
that the antipode $\un {S}$ of a Hopf algebra $B$ in a braided monoidal
category ${\mathcal C}$ is an antialgebra and anticoalgebra morphism, in the
sense that
\begin{equation}\label{antiac}
\un {S}\circ \un{m}=\un{m}\circ (\un{S}\ot \un{S})\circ c_{B, B}
\mbox{${\;\;}$and${\;\;}$}
\un{\Delta }\circ \un{S}=c_{B, B}\circ (\un{S}\ot \un{S})\circ \un{\Delta }.
\end{equation}
Recall also that an algebra $A$ in a braided monoidal 
category ${\mathcal C}$ is called
quantum commutative if $\un{m}\circ c_{A, A}=\un{m}$.\\

Assume that $(H, \Delta, \va , \Phi)$
is a quasi-bialgebra, and let $U, V, W$ be left $H$-modules. We define a
left $H$-action on $U\ot V$ by
$$h\cdot (u\ot v)=\sum h_1\cd u\ot h_2\cd v.$$
We have isomorphisms $a_{U, V, W}:\ (U\ot V)\ot W\ra
U\ot (V\ot W)$ in ${}_H{\mathcal M}$ given by
\begin{equation} \label{as}
a_{U, V, W}((u\ot v)\ot w)= \Phi \cd (u\ot (v\ot w)).
\end{equation}
The counit $\varepsilon:\ H\to k$ makes $k\in {}_H{\mathcal M}$, and
the natural isomorphisms $\lambda:\ k\ot H\to H$ and $\rho:\ H\ot
k\to H$ are in ${}_H{\mathcal M}$.
With this structures, $({}_H{\mathcal M}, \ot, k, a, \lambda, \rho)$ is
a monoidal category.\\
If $H$ is a quasi-Hopf algebra then the
category of finite dimensional left
$H$-modules is rigid; the left dual of $V$ is
$V^*$ with the $H$-module structure given by
$(h\cd \varphi )(v)=
\varphi (S(h)\cd v)$, for all $v\in V$,
$\varphi \in V^*$, $h\in H$ and with
\begin{eqnarray}
&&\hspace*{-2cm}
{\rm ev}_V(\varphi \ot v)=\varphi (\a \cd v),
\mbox{${\;\;\;}$}
{\rm coev}_V(1)=\sum _{i=1}^n\b \cd v_i\ot v^i,\label{qrig}
\end{eqnarray}
where $\{v_i\}$ is a basis in
$V$ with dual basis $\{v^i\}$.\\
Now let $H$ be a quasitriangular quasi-Hopf algebra,
with $R$-matrix $R=\sum R^1\ot R^2$. For two left $H$-modules $U$
and $V$, we define $$c_{U,V}:\ U\ot V\to V\ot U$$ by
\begin{equation}\label{br}
c_{U, V}(u\ot v)=\sum R^2\cd v\ot R^1\cd u
\end{equation}
and then $({}_H{\mathcal M}, \ot, k, a, \lambda, \rho, c)$ is a braided
monoidal category (cf. \cite{k} or \cite{m2}).

\section{Yetter-Drinfeld modules and the quasi-Yang-Baxter equation}\selabel{3}
From \cite{m1}, we recall the notion of Yetter-Drinfeld
module over a quasi-bialgebra.

\begin{definition}
Let $H$ be a quasi-bialgebra with reassociator $\Phi$. A
left $H$-module $M$ together with a left $H$-coaction
$$
\lambda_M:\ M\to H\ot M,~~\lambda_M(m)=\sum m_{(-1)}\ot m_{(0)}
$$
is called a left Yetter-Drinfeld module if the following
equalities hold, for all $h\in H$ and $m\in M$:
\begin{eqnarray}
&&\sum X^1m_{(-1)}\ot (X^2\cd m_{(0)})_{(-1)}X^3
\ot (X^2\cd m_{(0)})_{(0)}\nonumber\\%
&&\hspace*{1cm}=\sum X^1(Y^1\cd
m)_{(-1)_1}Y^2\ot X^2(Y^1\cd m)_{(-1)_2}Y^3\ot X^3\cd (Y^1\cd
m)_{(0)}\label{y1}\\%
&&\sum \va
(m_{(-1)})m_{(0)}=m\label{y2}\\
&&\sum h_1m_{(-1)}\ot h_2\cd m_{(0)}= \sum (h_1\cd
m)_{(-1)}h_2\ot (h_1\cd m)_{(0)}.\label{y3}
\end{eqnarray}
\end{definition}

The category of left Yetter-Drinfeld $H$-modules and $k$-linear
maps that intertwine the $H$-action and $H$-coaction is denoted by
$\yd$. In \cite{m1} it is shown that $\yd$ is a prebraided monoidal
category. The forgetful functor $\yd\to {}_H{\mathcal M}$ is monoidal,
and the coaction on the tensor product $M\ot N$ of two
Yetter-Drinfeld modules $M$ and $N$ is given by
\begin{eqnarray}
&&\hspace*{-2cm}
\lambda _{M\ot N}(m\ot n)=\sum X^1(x^1Y^1\cd m)_
{(-1)} x^2(Y^2\cd n)_{(-1)}Y^3\\
& \ot & X^2\cd (x^1Y^1\cd
m)_{(0)}\ot X^3x^3\cd (Y^2\cd n)_{(0)}.\label{y4}
\end{eqnarray}
The braiding is given by
\begin{equation}\label{y5}
c_{M, N}(m\ot n)=\sum m_{(-1)}\cd n\ot m_{(0)}.
\end{equation}
This braiding is invertible if $H$ is a quasi-Hopf algebra \cite{bn},
and its inverse is then given by
\begin{eqnarray}
&&\hspace*{-1cm}
c^{-1}_{M, N}(n\ot m)=\sum y^3_1X^2\cd (x^1\cd
m)_{(0)}\nonumber\\
&& \ot S^{-1}(S(y^1)\a y^2X^1(x^1\cd
m)_{(-1)}x^2\b S(y^3_2X^3x^3))\cd n.\label{y6}
\end{eqnarray}
Let $(H, R)$ be a quasitriangular quasi-bialgebra. It is well-known (see
for example \cite{k}) that $R$ satisfies the so-called
quasi-Yang-Baxter equation in $H\otimes H\otimes H$:
$$
R_{12}\Phi _{312}R_{13}\Phi ^{-1}_{132}R_{23}\Phi =
\Phi _{321}R_{23}\Phi ^{-1}_{231}R_{13}\Phi _{213}R_{12}.
$$
On the other hand, if $H$ is a bialgebra and $M$ is a left-right
Yetter-Drinfeld module over $H$, with structures
\begin{eqnarray*}
&&\hspace*{-2cm}
H\otimes M\rightarrow M,\;\;\;\;h\otimes m\mapsto h\cdot m;\\
&&\hspace*{-2cm}
M\rightarrow M\otimes H,\;\;\;\;m\mapsto \sum m_{(0)}\otimes m_{(1)},
\end{eqnarray*}
then the map $R_M:\ M\otimes M\rightarrow M\otimes M$, $R_M(m\otimes n)=
\sum n_{(1)}\cdot m\otimes n_{(0)}$ is a solution 
in $\End(M\otimes M\otimes M)$ of the quantum
Yang-Baxter equation
$$
R_{12}R_{13}R_{23}=R_{23}R_{13}R_{12},
$$
see for instance \cite{lr}.\\
We will show a similar result for quasi-bialgebras; 
 first we define left-right Yetter-Drinfeld modules over
quasi-bialgebras as follows
$$ {}_H{\mathcal YD}^H={} ^{H^{\rm cop}}_{H^{\rm cop}}{\mathcal YD}.$$ 
This is stated more explicitely in the next definition.

\begin{definition}
Let $H$ be a quasi-bialgebra. A $k$-linear space $M$ with a left
$H$-action $h\otimes m\mapsto h\cdot m$, and a right $H$-coaction
$M\rightarrow M\otimes H, \;m\mapsto
\sum m_{(0)}\otimes m_{(1)}$ is called a left-right Yetter-Drinfeld module
if the following relations hold, for all $m\in M$ and $h\in H$:
\begin{eqnarray}
&&
\sum (x^2\cdot m_{(0)})_{(0)}\otimes (x^2\cdot m_{(0)})_{(1)}x^1
\otimes x^3m_{(1)}\nonumber \\
&&\hspace*{1cm}=\sum x^1\cdot (y^3\cdot m)_{(0)}\otimes x^2(y^3\cdot m)_{(1)_1}y^1
\otimes x^3(y^3\cdot m)_{(1)_2}y^2\label{lry1}\\
&&\sum \varepsilon (m_{(1)})m_{(0)}=m\label{lry2}\\
&&
\sum h_1\cdot m_{(0)}\otimes h_2m_{(1)}=\sum (h_2\cdot m)_{(0)}
\otimes (h_2\cdot m)_{(1)}h_1\label{lry3}.
\end{eqnarray}
\end{definition}

\begin{proposition}
Let $H$ be a quasi-bialgebra and $M\in {}_H{\mathcal YD}^H$.
The map $R=R_M:\ M\otimes M\rightarrow
M\otimes M$, $R(m\otimes n)=\sum n_{(1)}\cdot m\otimes n_{(0)}$,
is a solution of the quasi-Yang-Baxter equation 
\begin{equation}\eqlabel{qqyb}
R_{12}\Phi _{312}R_{13}\Phi ^{-1}_{132}R_{23}\Phi =
\Phi _{321}R_{23}\Phi ^{-1}_{231}R_{13}\Phi _{213}R_{12}
\end{equation}
on $\End (M\otimes M\otimes M)$.
\end{proposition}

We considered $R_{12}, \Phi _{312}$, etc. as elements in
$\End (M\otimes M\otimes M)$ by left multiplication, for example
$R_{12}(l\otimes m\otimes n)=\sum R^1\cdot l\otimes R^2\cdot m\otimes n$,
$\Phi _{312}(l\otimes m\otimes n)=\sum X^2\cdot l\otimes X^3\cdot m\otimes
X^1\cdot n$ etc.

\begin{proof}
${}_H{\mathcal YD}^H$ is a prebraided category, hence
the result is a consequence of the fact (see \cite{k}) that the
braiding satisfies the categorical version of the Yang-Baxter equation.
A direct proof is also possible. For all $l, m, n\in M$, we compute
that
\begin{eqnarray*}
&&\hspace*{-2cm}
R_{12}\Phi _{312}R_{13}\Phi ^{-1}_{132}R_{23}\Phi (l\otimes m\otimes n)\\
&=&\sum (Y^3x^3(X^3\cdot n)_{(1)}X^2\cdot m)_{(1)}Y^2(x^2\cdot
(X^3\cdot n)_{(0)})_{(1)}x^1X^1\cdot l\\
&&~~~~\otimes (Y^3x^3(X^3\cdot n)_{(1)}X^2
\cdot m)_{(0)}\otimes Y^1\cdot (x^2\cdot (X^3\cdot n)_{(0)})_{(0)}\\
{\rm (\ref{lry1})}&=&
\sum (Y^3x^3(y^3X^3\cdot n)_{(1)_2}y^2X^2\cdot m)_{(1)}Y^2x^2
(y^3X^3\cdot n)_{(1)_1}y^1X^1\cdot l\\
&&~~~~\otimes (Y^3x^3(y^3X^3\cdot n)_{(1)_2}
y^2X^2\cdot m)_{(0)}\otimes Y^1x^1\cdot (y^3X^3\cdot n)_{(0)}\\
&=&\sum (n_{(1)_2}\cdot m)_{(1)}n_{(1)_1}\cdot l\otimes
(n_{(1)_2}\cdot m)_{(0)}\otimes n_{(0)}\\
{\rm (\ref{lry3})}&=&\sum n_{(1)_2}m_{(1)}\cdot l
\otimes n_{(1)_1}\cdot m_{(0)}\otimes n_{(0)}
\end{eqnarray*}
and
\begin{eqnarray*}
&&\hspace*{-2cm}
\Phi _{321}R_{23}\Phi ^{-1}_{231}R_{13}\Phi _{213}R_{12}
(l\otimes m\otimes n)\\
&=&\sum Y^3x^3(X^3\cdot n)_{(1)}X^2m_{(1)}\cdot l\otimes
Y^2(x^2\cdot (X^3\cdot n)_{(0)})_{(1)}x^1X^1\cdot m_{(0)}\\
&&~~~~\otimes
Y^1\cdot (x^2\cdot (X^3\cdot n)_{(0)})_{(0)}\\
{\rm (\ref{lry1})}&=&\sum Y^3x^3(y^3X^3\cdot n)_{(1)_2}y^2X^2m_{(1)}\cdot l\otimes
Y^2x^2(y^3X^3\cdot n)_{(1)_1}y^1X^1\cdot m_{(0)}\\
&&~~~~\otimes Y^1x^1\cdot (y^3X^3\cdot n)_{(0)}\\
&=&\sum n_{(1)_2}m_{(1)}\cdot l\otimes n_{(1)_1}\cdot m_{(0)}\otimes
n_{(0)}
\end{eqnarray*}
and \equref{qqyb} follows.
\end{proof}

We will now present a generalization of \cite[Prop. 4.4.2]{lr},
stating that the dual
$M^*$ of a finite dimensional right-left Yetter-Drinfeld module is a
left-right Yetter-Drinfeld module and that $R_{M^*}=R_M^*$.\\
First we define right-left Yetter-Drinfeld modules for
quasi-bialgebras as follows:
$${}^H{\mathcal YD}_H={} _{H^{\rm op, cop}}{\mathcal YD}^{H^{\rm op, cop}}.$$
More explicitely:

\begin{definition}
Let $H$ be a quasi-bialgebra. A $k$-linear space $M$ with a right
$H$-action $m\otimes h\mapsto m\cdot h$, and a left $H$-coaction
$M\rightarrow H\otimes M, \;m\mapsto
\sum m_{(-1)}\otimes m_{(0)}$ is called a right-left Yetter-Drinfeld module
if the following relations hold, for all $m\in M$ and $h\in H$:
\begin{eqnarray}
&&
\sum m_{(-1)}x^1\otimes x^3(m_{(0)}\cdot x^2)_{(-1)}\otimes
(m_{(0)}\cdot x^2)_{(0)}\nonumber\\
&&\hspace{1cm}=\sum y^2(m\cdot y^1)_{(-1)_1}x^1\otimes y^3(m\cdot y^1)_{(-1)_2}x^2
\otimes (m\cdot y^1)_{(0)}\cdot x^3\label{rly1}\\
&&
\sum \varepsilon (m_{(-1)})m_{(0)}=m\label{rly2}\\
&&
\sum m_{(-1)}h_1\otimes m_{(0)}\cdot h_2=
\sum h_2(m\cdot h_1)_{(-1)}\otimes (m\cdot h_1)_{(0)}\label{rly3}.
\end{eqnarray}
\end{definition}

For $M\in{}^H{\mathcal YD}_H$, we consider the map
$$R_M:\ M\otimes M
\rightarrow M\otimes M,~~~
R_M(m\otimes n)=\sum m\cdot n_{(-1)}\otimes n_{(0)}.$$
If we consider $M$ as an object in ${}_{H^{op, cop}}{\mathcal YD}^{H^{op, cop}}$,
then we obtain the same map $R_M$, so $R_M$ is also a solution of the
corresponding quasi-Yang-Baxter equation, which is obtained after replacing
$\Phi $ by $\Phi _{op, cop}=\Phi ^{321}$).\\
Now let $M$ be a finite dimensional right-left Yetter-Drinfeld module.
Then $M^*$ is a left $H$-module,
with action given by $(h\cdot m^*)(m)=m^*(m\cdot h)$, for all $h\in H, m\in M,
m^*\in M^*$. We also define a $k$-linear map $M^*\rightarrow M^*\otimes H$,
$m^*\mapsto \sum m^*_{(0)}\otimes m^*_{(1)}$, by the condition
\begin{equation}\label{mor}
\sum m^*_{(0)}(m)m^*_{(1)}=\sum m^*(m_{(0)})m_{(-1)}
\end{equation}
for all  $m\in M$. We can prove now the following result.

\begin{proposition}
Let $H$ be a quasi-bialgebra, $M$ a finite dimensional
right-left Yetter-Drinfeld module. Then
\begin{itemize}
\item [(i)] $M^*\in _H{\mathcal YD}^H$;
\item [(ii)] $R_{M^*}=R_M^*$.
\end{itemize}
\end{proposition}

\begin{proof}
$(i)$ We prove that (\ref{lry1}), (\ref{lry2}), (\ref{lry3}) are satisfied.
For $m^{*}\in M^*$ and $m\in M$, we compute:
\begin{eqnarray*}
&&\hspace*{-2cm}
\sum (x^2\cdot m^*_{(0)})_{(0)}(m)(x^2\cdot m^*_{(0)})_{(1)}x^1\otimes
x^3m^*_{(1)}\\
{\rm (\ref{mor})}
&=&\sum (x^2\cdot m^*_{(0)})(m_{(0)})m_{(-1)}x^1\otimes x^3m^*_{(1)}\\
&=&\sum m^*_{(0)}(m_{(0)}\cdot x^2)m_{(-1)}x^1\otimes x^3m^*_{(1)}\\
{\rm (\ref{rly1})}&=&\sum m^*((m\cdot y^1)_{(0)}\cdot x^3)y^2
(m\cdot y^1)_{(-1)_1}x^1\otimes
y^3(m\cdot y^1)_{(-1)_2}x^2\\
&=&\sum (x^3\cdot m^*)((m\cdot y^1)_{(0)})y^2(m\cdot y^1)_{(-1)_1}x^1
\otimes y^3(m\cdot y^1)_{(-1)_2}x^2\\
{\rm (\ref{mor})}&=&\sum (x^3\cdot m^*)_{(0)}(m\cdot y^1)y^2
(x^3\cdot m^*)_{(1)_1}x^1
\otimes y^3(x^3\cdot m^*)_{(1)_2}x^2\\
&=&\sum (y^1\cdot (x^3\cdot m^*)_{(0)})(m)y^2(x^3\cdot m^*)_{(1)_1}x^1
\otimes y^3(x^3\cdot m^*)_{(1)_2}x^2
\end{eqnarray*}
so  obtain (\ref{lry1}). Now we compute:
\begin{eqnarray*}
&&\hspace*{-2cm}
\sum \varepsilon (m^*_{(1)})m^*_{(0)}(m)
=\sum \varepsilon (m^*_{(0)}(m)m^*_{(1)})\\
{\rm (\ref{mor})}&=&\sum \varepsilon (m^*(m_{(0)})m_{(-1)})
=\sum m^*(\varepsilon (m_{(-1)})m_{(0)})=m^*(m),
\end{eqnarray*}
using (\ref{rly2}) at the last step. Thus  (\ref{lry2}) holds. For $h\in H$, we compute:
\begin{eqnarray*}
&&\hspace*{-2cm}
\sum (h_1\cdot m^*_{(0)})(m)h_2m^*_{(1)}
=\sum m^*_{(0)}(m\cdot h_1)h_2m^*_{(1)}\\
{\rm (\ref{mor})}&=&\sum m^*((m\cdot h_1)_{(0)})h_2(m\cdot h_1)_{(-1)}\\
{\rm (\ref{rly3})}&=&\sum m^*(m_{(0)}\cdot h_2)m_{(-1)}h_1\\
&=&\sum (h_2\cdot m^*)(m_{(0)})m_{(-1)}h_1\\
&=&\sum (h_2\cdot m^*)_{(0)}(m)(h_2\cdot m^*)_{(1)}h_1
\end{eqnarray*}
and (\ref{lry3}) follows.\\
$(ii)$ We identify $(M\otimes M)^*=M^*\otimes M^*$, and we prove that
$R_{M^*}$ and $R_M^*$ coincide as maps $M^*\otimes M^*\rightarrow
M^*\otimes M^*$. For $m, n\in M$ and $m^*, n^*\in M^*$, we compute:
\begin{eqnarray*}
&&\hspace*{-2cm}
R_{M^*}(m^*\otimes n^*)(m\otimes n)=
\sum (n^*_{(1)}\cdot m^*)(m)n^*_{(0)}(n)\\
&=&\sum m^*(m\cdot n^*_{(1)})n^*_{(0)}(n)\\
{\rm (\ref{mor})}&=&\sum m^*(m\cdot n_{(-1)})n^*(n_{(0)})\\
&=&(m^*\otimes n^*)(R_M(m\otimes n))\\
&=& R_M^*(m^*\otimes n^*)(m\otimes n),
\end{eqnarray*}
as needed.
\end{proof}

\section{The quantum commutativity of $H_0$}\selabel{4}
Let $H$ be a Hopf algebra. It is well-known that $H$ is an algebra
in the monoidal category $\yd $, with left action and coaction given by
$$
h\triangleright h'=\sum h_1h'S(h_2),
\mbox{${\;\;\;}$}
\lambda (h)=\sum h_1\otimes h_2.
$$
Moreover, $H$ is quantum commutative as an algebra
in $\yd $, see for example \cite{cvoz}.\\
We will now prove a similar result for quasi-Hopf algebras. Let $H$ be a
quasi-Hopf algebra. In \cite{bpv}, a new multiplication
on $H$ was introduced; this multiplication is given by the formula
\begin{equation}\label{ma}
h\circ h^{'}=\sum X^1hS(x^1X^2)\a
x^2X^3_1h^{'}S(x^3X^3_2)
\end{equation}
for all $h, h'\in H$. $\beta$ is a unit for this multiplication $\circ$.
Let $H_0$ be the $k$-linear space $H$, with multiplication
$\circ$, and left $H$-action given by
\begin{equation}\label{s1}
h\tr h^{'}=\sum h_1h^{'}S(h_2).
\end{equation}
Then $H_0$ is a left $H$-module algebra. In $H_0$, we also define a left $H$-coaction,
as follows
\begin{eqnarray}
&&\hspace*{-2cm}
\lambda _{H_0}(h)=\sum h_{(-1)}\ot h_{(0)}\nonumber\\
&=&
\sum X^1Y^1_1h_1g^1S(q^2Y^2_2)Y^3\ot X^2Y^1_2h_2g^2S(X^3q^1Y^2_1),\label{s2}
\end{eqnarray}
where $f^{-1}=\sum g^1\otimes g^2$ and $q_R=\sum q^1
\otimes q^2$ are the elements defined by (\ref{g}) and (\ref{qr}).
Then $H_0$ is an algebra in $\yd$, see \cite{bn} for details.
In \prref{4.2}, we will show that 
$H_0$ is quantum commutative. But first we need the following formulas,
which are of independent interest. Recall that
$q_R=\sum q^1\otimes q^2$, $q_L$, $f=\sum f^1\ot f^2$
and $f^{-1}=\sum g^1\ot g^2$ are defined by
(\ref{qra}), (\ref{qla}), (\ref{f}) and (\ref{g}).

\begin{lemma}\lelabel{4.1}
Let $H$ be a quasi-Hopf algebra. Then we have
\begin{eqnarray}
&&\sum q^1y^1\otimes S(q^2y^2)y^3=1\otimes \alpha ,\label{l1}\\
&&\Phi (\Delta \otimes id)(f^{-1})=\sum g^1S(X^3)f^1\otimes g^2_1G^1
S(X^2)f^2\otimes g^2_2G^2S(X^1),\label{l2}\\
&&\sum S(g^1)\alpha g^2=S(\beta ),
\mbox{${\;\;}$}
\sum f^1\b S(f^2)=S(\a ),\label{l3}\\
&&
\sum S(q^2_2X^3)f^1\otimes S(q^1X^1\b S(q^2_1X^2)f^2)
=(id\ot S)(q_L),\label{l5}
\end{eqnarray}
\end{lemma}
 
\begin{proof}
(\ref{l1}) and (\ref{l2}) are a direct consequence of (\ref{qr})
and (\ref{pf}). (\ref{l3}) has been proved
in \cite[Lemma 2.6]{bn2} and \cite[Lemma 2.5]{bn3}.
We are left to prove (\ref{l5}). Using (\ref{qr2}), we obtain:
$$
(id \otimes \Delta )(q)=\sum (1\otimes S^{-1}(x^3g^2)\otimes S^{-1}(x^2g^1))
(q\otimes 1)(\Delta \otimes id)(q)
\Phi ^{-1}(id \otimes \Delta )(\Delta (x^1))
$$
and, using the formula (see \cite{bc1})
$$
(\Delta \otimes id)(q)\Phi ^{-1}=\sum Y^1\otimes q^1Y^2_1\otimes
S^{-1}(Y^3)q^2Y^2_2,
$$
we obtain
\begin{equation}\label{fox}
(id \otimes \Delta )(q)=\sum Q^1Y^1x^1_1\otimes S^{-1}(x^3g^2)
Q^2q^1Y^2_1x^1_{(2,1)}\otimes
S^{-1}(Y^3x^2g^1)q^2Y^2_2x^1_{(2,2)}
\end{equation}
where $q_R=\sum q^1\otimes q^2=\sum Q^1\otimes Q^2$.
Now we compute
\begin{eqnarray*}
&&\hspace*{-2cm}
\sum S(q^2_2X^3)f^1\otimes S(q^1X^1\b S(q^2_1X^2)f^2)\\
{\rm (\ref{fox})}&=&\sum S(q^2Y^2_2x^1_{(2, 2)}X^3)Y^3x^2
\otimes S(Q^1Y^1x^1_1X^1\b S(Q^2q^1Y^2_1x^1_{(2, 1)}X^2)x^3)\\
{\rm (\ref{q1})}&=&\sum S(q^2Y^2_2X^3x^1_2)Y^3x^2
\otimes S(Q^1Y^1X^1x^1_{(1, 1)}\b S(Q^2q^1Y^2_1X^2x^1_{(1, 2)})x^3)\\
{\rm (\ref{q5})}&=&\sum S(q^2Y^2_2X^3x^1)Y^3x^2\ot
S(Q^1Y^1X^1\b S(Q^2q^1Y^2_1X^2)x^3)\\
{\rm (\ref{q3})}&=&\sum
S(q^2y^2X^3_1Y^2x^1)y^3X^3_2Y^3x^2\ot
S(Q^1X^1Y^1_1\b S(Q^2q^1y^1X^2Y^1_2)x^3)\\
{\rm (\ref{q5}, \ref{l1})}&=&\sum
S(X^3_1x^1)\a X^3_2x^2\ot S(Q^1X^1\b S(Q^2X^2)x^3)\\
{\rm (\ref{q5}, \ref{q7})}&=&\sum S(x^1)\a x^2\ot
S(Q^1\b S(Q^2)x^3)\\
{\rm (\ref{qr}, \ref{q6})}&=&\sum S(x^1)\a x^2\ot S(x^3),
\end{eqnarray*}
as needed.
\end{proof}

We can prove now the main result of this Section.

\begin{proposition}\prlabel{4.2}
Let $H$ be a quasi-Hopf algebra. Then $H_0$ is
quantum commutative as an algebra in $\yd $, that is,
for all $h, h^{'}\in H$:
$$
h\circ h'=\sum (h_{(-1)}\triangleright h')\circ h_{(0)}.
$$
\end{proposition}

\begin{proof}
For all $h, h^{'}\in H$ we compute:
\begin{eqnarray*}
&&\hspace*{-2cm}
\sum (h_{(-1)}\triangleright h')\circ h_{(0)}\\
{\rm (\ref{s2})}&=&\sum (X^1Y^1_1h_1g^1S(q^2Y^2_2)Y^3
\triangleright h')\circ
X^2Y^1_2h_2g^2S(X^3q^1Y^2_1)\\
{\rm (\ref{s1}, \ref{ma})}&=&\sum
Z^1X^1_1Y^1_{(1,1)}h_{(1,1)}g^1_1S(q^2Y^2_2)_1Y^3_1h'\\
&&\hspace*{15mm}S(x^1Z^2X^1_2Y^1_{(1,2)}h_{(1,2)}g^1_2S(q^2Y^2_2)_2Y^3_2)\\
&&\hspace*{15mm}\alpha x^2Z^3_1X^2Y^1_2h_2g^2S(x^3Z^3_2X^3q^1Y^2_1)\\
{\rm (\ref{q3}, \ref{q5})}&=&\sum
Z^1Y^1_{(1, 1)}h_{(1, 1)}g^1_1S(q^2Y^2_2)_1Y^3_1h^{'}\\
&&\hspace*{15mm}
S(Z^2Y^1_{(1, 2)}h_{(1, 2)}g^1_2S(q^2Y^2_2)_2Y^3_2)\\
&&\hspace*{15mm}\a Z^3Y^1_2h_2g^2S(q^1Y^2_1)\\
{\rm (\ref{ca})}&=&\sum Z^1[Y^1hS(Y^2)]_{(1, 1)}g^1_1S(q^2)_1
Y^3_1h^{'}\\
&&\hspace*{15mm}~S(Z^2[Y^1hS(Y^2)]_{(1, 2)}g^1_2S(q^2)_2Y^3_2)\\
&&\hspace*{15mm}\a Z^3[Y^1hS(Y^2)]_2g^2S(q^1)\\
{\rm (\ref{q1}, \ref{q5})}&=&\sum
Y^1hS(Y^2)Z^1g^1_1S(q^2)_1Y^3_1h^{'}S(Z^2g^1_2S(q^2)_2Y^3_2)
\a Z^3g^2S(q^1)\\
{\rm (\ref{l2})}&=&\sum Y^1hS(Y^2)g^1S(X^3)f^1S(q^2)_1Y^3_1h^{'}\\
&&\hspace*{15mm}S(g^2_1G^1S(X^2)f^2S(q^2)_2Y^3_2)\a g^2_2G^2S(q^1X^1)
\end{eqnarray*}
\begin{eqnarray*}
{\rm (\ref{q5}, \ref{l3})}&=&\sum
Y^1hS(X^3Y^2)f^1S(q^2)_1Y^3_1h^{'}
S(q^1X^1\b S(X^2)f^2S(q^2)_2Y^3_2)\\
{\rm (\ref{ca})}&=&\sum
Y^1hS(q^2_2X^3Y^2)f^1Y^3_1h^{'}
S(q^1X^1\b S(q^2_1X^2)f^2Y^3_2)\\
{\rm (\ref{l5})}&=&\sum Y^1hS(x^1Y^2)\a x^2Y^3_1h^{'}S(x^3Y^3_2)\\
{\rm (\ref{ma})}&=&h\circ h^{'}.
\end{eqnarray*}
\end{proof}

\section{Hopf modules in $\yd $. Integrals}\selabel{5}
Let $H$ be a quasi-Hopf algebra. The aim of this Section
is to define the space of integrals of a finite dimensional braided Hopf algebra
in $\yd$, and to prove, following \cite{ta}, \cite{doi}, that it is an object
of $\yd$, and that it has dimension $1$. We will apply our results to
the braided Hopf algebra associated to $H$, in the case where $H$ is a
quasitriangular quasi-Hopf algebra.\\
Let $A$ be an algebra in a monoidal category ${\mathcal C}$. Recall that
a right $A$-module $M$
is an object $M\in {\mathcal C}$ together with a morphism $\un {\omega }_M:\ M\ot A\ra M$
in ${\mathcal C}$ such that $\un {\omega }_M\circ (id_M\ot \un {\eta })=l_M^{-1}$
and the following diagram is commutative:
$$\begin{diagram}
(M\ot A)\ot A&\rTo^{{\un \omega }_M\ot id_A}&M\ot A&\rTo^{\un {\omega }_M}&M \\
\dTo^{a_{M, A, A}}&                         &      &      &\uTo_{\un {\omega }_M}\\
M\ot (A\ot A)&        &\rTo^{id_M\ot \un{m}}&                &M\ot A .
\end{diagram}$$
Clearly $A$ itself is a right $A$-module, by right multiplication.
Right comodules over a coalgebra $C$ in ${\mathcal C}$ can be defined in a similar
way: we need $N\in {\mathcal C}$
together with a morphism $\un {\rho }_N:\ N\ra N\ot C$ in ${\mathcal C}$ such that
$(id_N\ot \un {\va })\circ \un {\rho }_N=l_N$ and the following diagram
is commutative:
$$\begin{diagram}
N&\rTo^{{\un \rho }_N}&N\ot C&\rTo^{\un {\rho }_N\ot id_C}&(N\ot C)\ot C \\
\dTo^{\un {\rho }_N}&                         &      &      &\dTo_{a_{N, C, C}}\\
N\ot C&        &\rTo^{id_N\ot \un {\Delta }}&                &N\ot (C\ot C) .
\end{diagram}$$
$C$ itself is a right $C$-comodule via the
comultiplication $\un {\Delta }$.\\
From \cite{bkl}, \cite{m3}, \cite{ta}, we recall the following.

\begin{definition}
Let $B$ be a bialgebra in a braided category ${\mathcal C}$.
A right $B$-Hopf module is a triple $(M, \un {\omega }_M, \un {\r }_M)$, where
$(M, \un {\omega }_M)$ is a right $B$-module and $(M, \un {\r }_M)$ is a right $B$-comodule
such that $\un {\r }_M:\ M\ra M\ot B$ is right $B$-linear. The $B$-module structure
$\un {\omega }_{M\ot B}:\ (M\ot B)\ot B\ra M\ot B$ on
$M\ot B$ is given by the following composition:
\begin{equation}\label{hmc}
\begin{array}{ccc}
(M\ot B)\ot B&\rTo^{id_{M\ot B}\ot {\un \Delta }}&
(M\ot B)\ot (B\ot B)\\
&\rTo^{a_{M, B, B\ot B}}&M\ot (B\ot (B\ot B))\\
&\rTo^{id_M\ot a^{-1}_{B, B, B}}&M\ot ((B\ot B)\ot B)\\
&\rTo^{id_M\ot (c_{B, B}\ot id_B)}&M\ot ((B\ot B)\ot B)\\
&\rTo^{id_M\ot a_{B, B, B}}&M\ot (B\ot (B\ot B))\\
&\rTo^{a^{-1}_{M, B, B\ot B}}&(M\ot B)\ot (B\ot B)\\
&\rTo^{{\un \omega }_M\ot \un {m}}&M\ot B
\end{array}
\end{equation}
${\mathcal M}_B^B$ will denote the category of right $B$-Hopf modules
and  morphisms in ${\mathcal C}$ preserving the $B$-action and
the corresponding $B$-coaction.
\end{definition}

We can consider algebras, coalgebras, bialgebras and Hopf algebras in
the braided category $\yd$ over a quasi-Hopf algebra $H$. More precisely,
an algebra $B$ in $\yd$ is an object $B\in \yd$ such that
\begin{enumerate}
\item[-] $B$ is a left $H$-module algebra, i.e.
$B$ has a multiplication $\un {m}$
and a usual unit $1_B$ satisfying the
following conditions:
\begin{equation}\label{mal}
(ab)c=\sum (X^1\cd a)[(X^2\cd b)(X^3\cd c)],
\end{equation}
\begin{equation}
h\cd (ab)=\sum (h_1\cd a)(h_2\cd b),
\mbox{${\;\;}$}
h\cd 1_B=\va (h)1_B,
\end{equation}
for all $a, b, c\in B$ and $h\in H$.
\item[-] $B$ is a quasi-comodule
algebra, that is, the multiplication $\un {m}$ and the unit
$\un {\eta }$ of $B$ intertwine the $H$-coaction $\lambda _B$. By (\ref{y4})
this means:
\begin{eqnarray}
&&\hspace*{-2cm}
\lambda _B(bb^{'})=\sum X^1(x^1Y^1\cd b)_{(-1)}x^2(Y^2\cd b^{'})_{(-1)}Y^3 \nonumber\\
&&\ot [X^2\cd (x^1Y^1\cd b)_{(0)}]
[X^3x^3\cd (Y^2\cd b^{'})_{(0)}],\label{qca1}
\end{eqnarray}
for all $b, b^{'}\in B$, and
\begin{equation}\label{qca2}
\lambda _B(1_B)=1_H\ot 1_B.
\end{equation}
\end{enumerate}
$M\in \yd $ is a right $B$-module
if there exists a morphism $\un {\omega }_M:\ M\ot B\ra M$ in $\yd $ (we will denote
$\un {\omega }_M(m\ot b):=m\leftarrow b$) such that
\begin{equation}\label{rm1}
m\leftarrow 1_B=m,
\mbox{${\;\;}$}
(m\leftarrow b)\leftarrow b^{'}=\sum (X^1\cd m)\leftarrow [(X^2\cd b)(X^3\cd b^{'})]
\end{equation}
for all $m\in M$, $b, b^{'}\in B$.
The fact that $\un {\omega }_M$ is a morphism in $\yd $ means (see (\ref{y4}))
\begin{eqnarray}
&&\hspace*{-2cm}h\cd (m\leftarrow b)=\sum (h_1\cd m)\leftarrow (h_2\cd b),\label{olhl}\\
&&\hspace*{-2cm}\l _M(m\leftarrow b)=\sum X^1(x^1Y^1\cd m)_{(-1)}x^2(Y^2\cd b)_{(-1)}Y^3\nonumber\\
&&\ot [X^2\cd (x^1Y^1\cd m)_{(0)}]\leftarrow [X^3x^3\cd (Y^2\cd b)_{(0)}]\label{olhcol}
\end{eqnarray}
for all $m\in M$, $b\in B$.\\

Similarly,  $B\in \yd$ is a coalgebra if 
\begin{enumerate}
\item[-] $B$ is a left $H$-module coalgebra, i.e. $B$ has a comultiplication
$\un {\Delta }_B :\ B\ra B\ot B$ (we will denote
$\un {\Delta }(b)=\sum b_{\un {1}}\ot b_{\un {2}}$)
and a usual counit $\un {\va }_B$ such that:
\begin{eqnarray}
&&\hspace*{-2cm}\sum X^1\cd b_{(\un{1}, \un{1})}
\ot X^2\cd b_{(\un{1}, \un{2})}
\ot X^3\cd b_{\un{2}}=\sum b_{\un{1}}\ot b_{(\un{2}, \un{1})}\ot
b_{(\un{2}, \un{2})},\label{mc1}\\
&&\hspace*{-2cm}
\un {\Delta }_B(h\cd b)=\sum h_1\cd b_{\un{1}}\ot h_2\cd b_{\un{2}},
\mbox{${\;\;\;}$}
\un{\va }_B(h\cd b)=\va (h)\un{\va }_B(b),\label{mc2}
\end{eqnarray}
for all $h\in H$, $b\in B$, where we use
  the same notation for the quasi-coassociativity
of $\un {\Delta }_B$ as in \seref{2}.
\item[-] $B$ is a quasi-comodule coalgebra,
i.e. the comultiplication
$\un {\Delta }_B$ and the counit $\un {\va }_B$ intertwine the
$H$-coaction $\l _B$. Explicitly, for all $b\in B$ we must have that:
\begin{eqnarray}
&&\hspace*{-1cm}\sum b_{(-1)}\ot b_{(0)_{\un{1}}}\ot b_{(0)_{\un{2}}}
=\sum X^1(x^1Y^1\cd b_{\un{1}})_{(-1)}
x^2(Y^2\cd b_{\un{2}})_{(-1)}Y^3\nonumber\\
&&\ot
X^2\cd (x^1Y^1\cd b_{\un{1}})_{(0)}\ot
X^3x^3\cd (Y^2\cd b_{\un{2}})_{(0)},\label{qcc1}
\end{eqnarray}
and
\begin{equation}
\sum \va _B(b_{(0)})b_{(-1)}=\va _B(b)1.\label{qcc2}
\end{equation}
\end{enumerate}
A right $B$-comodule in $\yd $ is an object $M\in \yd$
together with a morphism $\un {\r }_{M}:\ M\ra M\ot B$ in $\yd $ (we will
denote $\un {\r }_M(m)=\sum m_{(\un {0})}\ot m_{(\un {1})}$ for all $m\in M$) such that
 the following relations hold, for all $m\in M$:
\begin{equation}
\sum X^1\cd m_{({\un 0},{\un 0})}\ot X^2\cd m_{({\un 0}, {\un 1})}\ot X^3\cd m_{({\un 1})}=
\sum m_{({\un 0})}\ot m_{{({\un 1})}_{\un 1}}\ot m_{{({\un 1})}_{\un 2}},\label{rbc1}
\end{equation}
\begin{equation}
\sum \un{\va }(m_{({\un 1})})m_{({\un 0})}=m,\label{rbc2}
\end{equation}
where we will denote
$$
(\un {\r }_M\ot id_B)(\un {\r }_M(m))=\sum m_{({\un 0}, {\un 0})}\ot
m_{({\un 0}, {\un 1})}\ot m_{({\un 1})}
\mbox{${\;\;}$etc. }
$$
The fact that $\un {\r }_M$ is a morphism in $\yd $ means that (see (\ref{y4}))
\begin{equation}
\un {\r }_M(h\cd m)=\sum h_1\cd m_{({\un 0})}\ot h_2\cd m_{({\un 1})},\label{rhl}
\end{equation}
and
\begin{eqnarray}
&&\hspace*{-2cm}\sum m_{(-1)}\ot m_{(0)_{({\un 0})}}\ot m_{(0)_{({\un 1})}}=\sum
X^1(x^1Y^1\cd m_{({\un 0})})_{(-1)}x^2(Y^2\cd m_{({\un 1})})_{(-1)}Y^3\nonumber \\
&&\ot X^2\cd (x^1Y^1\cd m_{({\un 0})})_{(0)}\ot X^3x^3\cd (Y^2\cd m_{({\un 1})})_{(0)},
\label{rlhcol}
\end{eqnarray}
for all $h\in H$ and $m\in M$.\\
Now, a bialgebra $B\in \yd$ is an algebra and a coalgebra in $\yd $ such that
$\un {\Delta }_B$ is an algebra morphism, i.e.
$\un {\Delta }_B(1_B)=1_B\ot 1_B$ and, by (\ref{bi})
and (\ref{y5}), for all $b, b^{'}\in B$ we have that:
\begin{eqnarray}
&&\hspace*{-2cm}
\Delta _B(bb^{'})=\sum [y^1X^1\cd b_{\un{1}}]
[y^2Y^1(x^1X^2\cd b_{\un{2}})_{(-1)}x^2
X^3_1\cd b^{'}_{\un{1}}]\nonumber\\
&&\ot [y^3_1Y^2\cd (x^1X^2\cd b_{\un{2}})_{(0)}]
[y^3_2Y^3x^3X^3_2\cd b^{'}_{\un{2}}].\label{by}
\end{eqnarray}
If $B\in \yd $ is a bialgebra then $M\in \yd $ is a right $B$-Hopf module
if $M$ is a right $B$-module (as above, we will denote
$\un {\omega }_M(m\ot b)=m\leftarrow b$) and a right $B$-comodule such that
the right $B$-coaction on $M$, $\un {\r }_M:\  M\ra M\ot B$, is right $B$-linear,
which means that the following relation holds, for all $m\in M$ and $b\in B$
(see (\ref{hmc})):
\begin{eqnarray}
&&\hspace*{-2cm}\un {\r }_M(m\leftarrow b)=\sum (y^1X^1\cd m_{({\un 0})})\leftarrow
[y^2Y^1(x^1X^2\cd m_{({\un 1})})_{(-1)}x^2X^3_1\cd b_{\un 1}]\nonumber\\
&&\ot [y^3_1Y^2\cd (x^1X^2\cd m_{({\un 1})})_{(0)}][y^3_2Y^3x^3X^3_2\cd b_{\un 2}].\label{hmyd}
\end{eqnarray}
Finally, a bialgebra $B$ in $\yd $ is a braided Hopf algebra
if there exists a morphism $\un {S}:\ B\ra B$ in $\yd $ such that
$\sum \un{S}(b_{\un{1}})b_{\un{2}}=\sum b_{\un{1}}S(b_{\un{2}})=
\un{\va }(b)1_B$, for all $b\in B$. Since $\un {S}$ is a morphism in $\yd $, we have that
\begin{equation}\label{smorf}
\un {S}(h\cd b)=h\cd \un {S}(b)
\mbox{${\;\;}$and${\;\;}$}
\sum \un{S}(b)_{(-1)}\ot \un{S}(b)_{(0)}=\sum b_{(-1)}\ot \un{S}(b_{(0)}),
\end{equation}
for all  $h\in H$, $b\in B$. Also, by (\ref{antiac}) and (\ref{y5}) we obtain that
\begin{equation}\label{santi}
\un{S}(bb^{'})=\sum [b_{(-1)}\cd \un{S}(b^{'})]{\un S}(b_{(0)})
\mbox{${\;\;}$and${\;\;}$}
{\un \Delta }({\un S}(b))=\sum b_{{\un 1}_{(-1)}}\cd {\un S}(b_{\un 2})\ot
{\un S}(b_{{\un 1}_{(0)}}),
\end{equation}
for all $b, b^{'}\in B$.\\
The first step to prove the existence and uniqueness of integrals in a finite
dimensional braided Hopf algebra is the structure theorem for Hopf modules. To
this end we need first the following result.

\begin{lemma}\lelabel{5.2}
Let $H$ be a quasi-bialgebra, $B$ a bialgebra in $\yd $ and $N\in \yd $. Then
$N\ot B\in {\mathcal M}^B_B$ with following action
$\un {\omega }_{N\ot B}:\ (N\ot B)\ot B\ra N\ot B$
 and coaction $\un{\r }_{N\ot B}:\ N\ot B\ra (N\ot B)\ot B$ given by
\begin{eqnarray}
&&
(n\ot b)\prec b^{'}=\sum X^1\cd n\ot
[(X^2\cd b)(X^3\cd b^{'})],\label{fst1}\\
&&
\un{\r }_{N\ot B}(n\ot b):=\sum x^1\cd n\ot x^2\cd b_{\un 1}\ot x^3\cd b_{\un 2},\label{fst2}
\end{eqnarray}
for all $n\in N$ and $b, b^{'}\in B$.
\end{lemma}

\begin{proof}
$\yd $ is a braided category, so $N\ot B\in \yd $.
It is not hard to see that (\ref{q1}) and (\ref{mal}) imply that
$\un{\omega }_{N\ot B}$ is left $H$-linear. It intertwines also the corresponding $H$-coaction.
Indeed, by (\ref{y4}), the left $H$-coaction on $(N\ot B)\ot B$ is given by
\begin{eqnarray*}
&&\hspace*{-8mm}
\l _{(N\ot B)\ot B}((n\ot b)\ot b^{'})\\
&=&\sum Z^1X^1(x^1Y^1y^1_1T^1_1\cd n)_{(-1)}x^2(Y^2y^1_2T^1_2\cd b)_{(-1)}Y^3y^2
(T^2\cd b^{'})_{(-1)}T^3\\
&&\hspace*{5mm} \ot Z^2_1X^2\cd (x^1Y^1y^1_1T^1_1\cd n)_{(0)}\ot
Z^2_2X^3x^3\cd (Y^2y^1_2T^1_2\cd b)_{(0)}\ot Z^3y^3\cd (T^2\cd b^{'})_{(0)},
\end{eqnarray*}
for all $n\in N$, $b, b^{'}\in B$. Therefore:
\begin{eqnarray*}
&&\hspace*{-2cm}
(id_H\ot \un{\omega }_{N\ot B})\circ \l _{(N\ot B)\ot B}((n\ot b)\ot b^{'})\\
{\rm (\ref{fst1})}&=&
\sum Z^1X^1(x^1Y^1y^1_1T^1_1\cd n)_{(-1)}x^2(Y^2y^1_2T^1_2\cd b)_{(-1)}\\
&&\hspace*{1cm}Y^3y^2(T^2\cd b^{'})_{(-1)}T^3
\ot W^1Z^2_1X^2\cd (x^1Y^1y^1_1T^1_1\cd n)_{(0)}\\
&&\hspace*{1cm}\ot
[(W^2Z^2_2X^3x^3\cd (Y^2y^1_2T^1_2\cd b)_{(0)}]
[(W^3Z^3y^3\cd (T^2\cd b^{'})_{(0)}]\\
{\rm (\ref{q3},\ref{y3},\ref{mal})}&=&\sum
Z^1(X^1_1x^1Y^1y^1_1T^1_1\cd n)_{(-1)}X^1_2x^2(Y^2y^1_2T^1_2\cd b)_{(-1)}\\
&&\hspace*{1cm}Y^3y^2(T^2\cd b^{'})_{(-1)}T^3
\ot Z^2\cd (X^1_1x^1Y^1y^1_1T^1_1\cd n)_{(0)}\\
&&\hspace*{1cm}\ot Z^3\cd
[(X^2x^3\cd (Y^2y^1_2T^1_2\cd b)_{(0)})
(X^3y^3\cd (T^2\cd b^{'})_{(0)})]\\
{\rm (\ref{q3})\mbox{${\;}$}twice, (\ref{y3})}&=&\sum
Z^1(x^1Y^1T^1_1\cd n)_{(-1)}x^2X^1(y^1Y^2T^1_2\cd b)_{(-1)}y^2\\
&&\hspace*{1cm}(Y^3_1T^2\cd b^{'})_{(-1)}Y^3_2T^3
\ot Z^2\cd (x^1Y^1T^1_1\cd n)_{(0)}\\
&&\hspace*{1cm}\ot
Z^3x^3\cd [(X^2\cd (y^1Y^2T^1_2\cd b)_{(0)})
(X^3y^3\cd (Y^3_1T^2\cd b^{'})_{(0)})]\\
{\rm (\ref{q3},\ref{qca1})}&=&\sum
Z^1(x^1Y^1T^1\cd n)_{(-1)}x^2[(Y^2_1T^2\cd b)(Y^2_2T^3\cd b^{'})]_{(-1)}Y^3\\
&&\hspace*{1cm}\ot Z^2\cd (x^1Y^1T^1\cd n)_{(0)}\ot Z^3x^3\cd
[(Y^2_1T^2\cd b)(Y^2_2T^3\cd b^{'})]_{(0)}\\
{\rm (\ref{y4},\ref{fst1})}&=&\sum
\l _{N\ot B}(T^1\cd n\ot (T^2\cd b)(T^3\cd b^{'}))\\
&=&\l _{N\ot B}\circ \un{\omega }_{N\ot B}((n\ot b)\ot b^{'})
\end{eqnarray*}
for all $n\in N$ and $b, b^{'}\in B$.
In a similar way, it can be proved that the map $\un{\r}_{N\ot B}$ is a morphism in $\yd$, we leave
it to the reader to verify the details.\\
Using (\ref{mal}) and (\ref{q3}), it easily follows that $N\ot B$ is a right $B$-module.
Also, it is not hard to see that (\ref{mc1}), (\ref{mc2})
and (\ref{q3}) imply that $N\ot B$ is a right $B$-comodule. It remains only to
show that ${\un {\r }}_{N\ot B}$ is right $B$-linear. By (\ref{hmc}), we
have that the right $B$-module structure of $(N\ot B)\ot B$ is given by
\begin{eqnarray*}
&&\hspace*{-15mm}
[(n\ot b)\ot b^{'}]\bullet b^{''}\\
&=&\sum [Z^1y^1_1X^1_1\cd n\ot (Z^2y^1_2X^1_2\cd b)(Z^3y^2Y^1(x^1X^2\cd b^{'})_{(-1)}
x^2X^3_1\cd b^{''}_{\un 1})]\\
&&\hspace*{1cm}\ot [y^3_1Y^2\cd (x^1X^2\cd b^{'})_{(0)}][y^3_2Y^3x^3X^3_2\cd b^{''}_{\un 2}],
\end{eqnarray*}
for all $n\in N$ and $b, b^{'}, b^{''}\in B$. This allows us to compute, for any $n\in N$ and
$b, b^{'}\in B$, that:
\begin{eqnarray*}
&&\hspace*{-2cm}
\un {\r }_{N\ot B}(n\ot b)\bullet b^{'}
=\sum [(z^1\cd n\ot z^2\cd b_{\un 1})\ot z^3\cd b_{\un 2}]\bullet b^{'}\\
&=&\sum [Z^1y^1_1X^1_1z^1\cd n\\
&&\hspace*{1cm} \ot(Z^2y^1_2X^1_2z^2\cd b_{\un 1})(Z^3y^2Y^1
(x^1X^2z^3\cd b_{\un 2})_{(-1)}x^2X^3_1\cd b^{'}_{\un 1})]\\
&&\hspace*{1cm}\ot [y^3_1Y^2\cd (x^1X^2z^3\cd b_{\un 2})_{(0)}][y^3_2Y^3x^3X^3_2\cd b^{'}_{\un 2}]
\end{eqnarray*}
\begin{eqnarray*}
{\rm (\ref{q3})}&=&\sum [Z^1y^1_1z^1X^1\cd n\ot (Z^2y^1_2z^2T^1X^2_1\cd b_{\un 1})\\
&&\hspace*{1cm}(Z^3y^2Y^1
(x^1z^3_1T^2X^2_2\cd b_{\un 2})_{(-1)}
x^2z^3_{(2, 1)}T^3_1X^3_1\cd b^{'}_{\un 1})]\\
&&\hspace*{1cm}\ot
[y^3_1Y^2\cd (x^1z^3_1T^2X^2_2\cd b_{\un 2})_{(0)}]
[y^3_2Y^3x^3z^3_{(2, 2)}T^3_2X^3_2\cd b^{'}_{\un 2}]\\
{\rm (\ref{q1},\ref{y3})}&=&\sum [Z^1y^1_1z^1X^1\cd n
\ot (Z^2y^1_2z^2T^1X^2_1\cd b_{\un 1})\\
&&\hspace*{1cm}
(Z^3y^2Y^1z^3_{(1, 1)}(x^1T^2X^2_2\cd b_{\un 2})_{(-1)}
x^2T^3_1X^3_1\cd b^{'}_{\un 1})]\\
&&\hspace*{1cm}\ot
[y^3_1Y^2z^3_{(1, 2)}\cd (x^1T^2X^2_2\cd b_{\un 2})_{(0)}]
[y^3_2Y^3z^3_2x^3T^3_2X^3_2\cd b^{'}_{\un 2}]\\
{\rm (\ref{q1},\ref{q3},\ref{mal})}&=&
\sum \{y^1X^1\cd n
\ot y^2\cd [(z^1T^1X^2_1\cd b_{\un 1})\\
&&\hspace*{1cm}
(z^2Y^1(x^1T^2X^2_2\cd b_{\un 2})_{(-1)}x^2T^3_1X^3_1\cd b^{'}_{\un 1})]\}\\
&&\hspace*{1cm}\ot y^3\cd \{[z^3_1Y^2\cd (x^1T^2X^2_2\cd b_{\un 2})_{(0)}]
[z^3_2Y^3x^3T^3_2X^3_2\cd b^{'}_{\un 2}]\}\\
{\rm (\ref{mc2},\ref{by})}&=&\sum \{y^1X^1\cd n\ot
y^2\cd [(X^2\cd b)(X^3\cd b^{'})]_{\un 1}\}\ot y^3\cd [(X^2\cd b)(X^3\cd b^{'})]_{\un 2}\\
{\rm (\ref{fst2},\ref{fst1})}&=&\sum
\un {\r }_{N\ot B}(X^1\cd n\ot (X^2\cd b)(X^3\cd b^{'}))
=\un {\r }_{N\ot B}((n\ot b)\prec b^{'}),
\end{eqnarray*}
as needed.
\end{proof}

Our next result is the Fundamental Theoreom for Hopf modules in the braided
monoidal category $\yd$, generalizing \cite[Theorem 1]{doi}.

\begin{theorem}\thlabel{4.3}
Let $H$ be a quasi-Hopf algebra, $B$ a Hopf algebra in $\yd $ and
$M\in {\mathcal M}_B^B$.
\begin{itemize}
\item[(i)] $M^{{\rm co}B}=\{m\in M\mid \un {\r }_M(m)=m\ot 1_B\}\in \yd $.
\item[(ii)] For all $m\in M$, we have that $P(m)=\sum m_{({\un 0})}\leftarrow
\un{S}(m_{({\un 1})})\in M^{{\rm co}B}$.
\item[(iii)] $\un{\r }_M(n\leftarrow b)=\sum (x^1\cd n)\leftarrow (x^2\cd b_{\un 1})
\ot x^3\cd b_{\un 2}$ and $P(n\leftarrow b)=\un{\va }(b)n$, for all $n\in M^{{\rm co}B}$ and
$b\in B$.
\item[(iv)] The map 
$$F:\ M^{{\rm co}B}\ot B\ra M,~~~F(n\ot b)=n\leftarrow b,$$
 is an isomorphism
of Hopf modules in $\yd $, with inverse $G$ given by
$$G(m)=\sum P(m_{({\un 0})})\ot m_{({\un 1})}.$$
\end{itemize}
\end{theorem}

\begin{proof}
(i)  If $n\in M^{{\rm co}B}$, then $\un{\r }_M(h\cd n)=\sum h_1\cd n\ot h_2\cd 1_B=
h\cd n\ot 1_B$, by (\ref{olhl}) and (\ref{mal}). This shows that $M^{{\rm co}B}$
is an $H$-submodule of $M$. On the other hand, for any $n\in N$ we have
\begin{eqnarray*}
&&\hspace*{-2cm}
\sum n_{(-1)}\ot n_{(0)_{({\un 0})}}\ot n_{(0)_{({\un 1})}}\\
{\rm (\ref{rlhcol})}&=&\sum X^1(x^1Y^1\cd n)_{(-1)}x^2(Y^2\cd 1_B)_{(-1)}Y^3\\
&&\hspace*{1cm} \ot X^2\cd (x^1Y^1\cd n)_{(0)}
\ot X^3x^3\cd (Y^2\cd 1_B)_{(0)}\\
{\rm (\ref{mal})\mbox{${\;}$twice,${\;}$}(\ref{qca2})}&=&\sum n_{(-1)}\ot n_{(0)}\ot 1_B.
\end{eqnarray*}
Thus, $\r _M(n)=\sum n_{(-1)}\ot n_{(0)}\in H\ot M^{{\rm co}B}$ which means that
$M^{{\rm co}B}$ is a left $H$-quasi-subcomodule of $M$. It follows from the above arguments
that $M^{{\rm co}B}\in \yd $.\\

(ii) For any $m\in M$, we have that
\begin{eqnarray*}
\un {\r }_M(P(m))
&=&\sum \un {\r }_M(m_{({\un 0})}\leftarrow \un{S}(m_{({\un 1})}))\\
{\rm (\ref{hmyd})}&=&\sum (y^1X^1\cd m_{({\un 0}, {\un 0})})\leftarrow
[y^2Y^1(x^1X^2\cd m_{({\un 0}, {\un 1})})_{(-1)}x^2X^3_1\cd
\un{S}(m_{({\un 1})})_{\un 1}]\\
&&\hspace*{5mm}[y^3_1Y^2\cd (x^1X^2\cd m_{({\un 0},{\un 1})})_{(0)}]
[y^3_2Y^3x^3X^3_2\cd \un{S}(m_{({\un 1})})_{\un 2}]\\
{\rm (\ref{rbc1},\ref{smorf},\ref{santi})}&=&\sum
(y^1\cd m_{({\un 0})})\leftarrow\\
&&\hspace*{5mm} [y^2Y^1(x^1\cd m_{({\un 1})_{\un 1}})_{(-1)}x^2
\un{S}(m_{({\un 1})_{({\un 2}, {\un 1})}})_{(-1)}\cd
\un{S}(m_{({\un 1})_{({\un 2}, {\un 2})}})]\\
&&\hspace*{5mm}\ot y^3\cd \{[Y^2\cd (x^1\cd m_{({\un 1})_{\un 1}})_{(0)}][Y^3x^3\cd
\un{S}(m_{({\un 1})_{({\un 2}, {\un 1})}})_{(0)}]\}\\
{\rm (\ref{qca1})}&=&\sum (y^1\cd m_{({\un 0})})\leftarrow\\
&&\hspace*{5mm}
y^2[(x^1\cd m_{({\un 1})_{\un 1}})(x^2\cd \un{S}(m_{({\un 1})_{({\un 2}, {\un 1})}}))]_{(-1)}
x^3\cd \un{S}(m_{({\un 1})_{({\un 2}, {\un 2})}})\\
&&\hspace*{5mm}\ot y^3\cd [(x^1\cd m_{({\un 1})_{\un 1}})
(x^2\cd \un{S}(m_{({\un 1})_{({\un 2}, {\un 1})}}))]_{(0)}\\
{\rm (\ref{smorf},\ref{mc1},\ref{mal})}&=&\sum (y^1\cd m_{({\un 0})})\leftarrow
(y^2\cd \un{S}(m_{({\un 1})}))\ot y^3\cd 1_B=P(m)\ot 1_B.
\end{eqnarray*}

(iii) For all $n\in N$ and $b\in B$, we compute, using (\ref{hmyd}),
\begin{eqnarray*}
\un{\r }_M(n\leftarrow b)&=&
\sum
(y^1X^1\cd n)\leftarrow [y^2Y^1(x^1X^2\cd 1_B)_{(-1)}x^2X^3_1\cd b_{\un 1}]\\
&&\hspace*{1cm}
\ot [y^3_1Y^2\cd (x^1X^2\cd 1_B)_{(0)}][y^3_2Y^3x^3X^3_2\cd b_{\un 2}]\\
{\rm (\ref{mal},\ref{qca2})}&=&\sum (y^1\cd n)\leftarrow
(y^2\cd b_{\un 1})\ot y^3\cd b_{\un 2}.
\end{eqnarray*}
For all $n\in M^{{\rm co}B}$, we find
\begin{eqnarray*}
P(n\leftarrow b)&=&\sum
[(y^1\cd n)\leftarrow (y^2\cd b_{\un 1})]\leftarrow \un{S}(y^3\cd b_{\un 2})\\
{\rm (\ref{rm1},\ref{smorf})}&=&\sum n\leftarrow b_{\un 1}\un{S}(b_{\un 2})
=\un{\va }(b)n\leftarrow 1_B
=\un{\va}(b)n.
\end{eqnarray*}

(iv) By (i) and \leref{5.2}, we obtain that $M^{{\rm co}B}\ot B\in {\mathcal M}_B^B$.
It follows from (\ref{olhl}) that $F$ is left $H$-linear. It also intertwines
the corresponding left $H$-coaction by (\ref{y4}) and (\ref{olhcol}). Now we will
prove that $F$ and $G$ are inverses. For all $m\in M$, we have
\begin{eqnarray*}
FG(m)&=&\sum P(m_{({\un 0})})\leftarrow m_{({\un 1})}\\
{\rm (\ref{rm1})}&=&\sum (X^1\cd m_{({\un 0}, {\un 0})})\leftarrow
[(X^2\cd \un{S}(m_{({\un 0}, {\un 1})}))(X^3\cd m_{({\un 1})})]\\
{\rm (\ref{smorf},\ref{rbc1},\ref{rbc2})}&=&\sum m_{({\un 0})}\leftarrow
\un{S}(m_{({\un 1})_{\un 1}})m_{({\un 1})_{\un 2}}=m\leftarrow 1_B=m.
\end{eqnarray*}
Similarly, for any $n\in M^{{\rm co}B}$ and $b\in B$, we compute
\begin{eqnarray*}
GF(n\ot b)&=&\sum P((n\leftarrow b)_{({\un 0})})\ot (n\leftarrow b)_{({\un 1})}\\
{\rm (iii)}&=&\sum P((x^1\cd n)\leftarrow (x^2\cd b_{\un 1}))\ot x^3\cd b_{\un 2}\\
{\rm (iii), (\ref{mc2})}&=&\sum P(n)\ot b=n\ot b.
\end{eqnarray*}
We are left to show that $F$ is a morphism in ${\mathcal M}_B^B$. It is not hard to see
that (\ref{fst1}) and (\ref{rm1}) imply that $F$ is right $B$-linear. Also, (iii) implies
that
$$
\un{\r }_M\circ F(n\ot b)=(F\ot id_B)\circ \un{\r }_{M^{{\rm co}B}\ot B}(n\ot b)=\sum
(x^1\cd n)\leftarrow (x^2\cd b_{\un 1})\ot x^3\cd b_{\un 2},
$$
for all $n\in N$ and $b\in B$, and this finishes the proof.
\end{proof}

Let $H$ be a quasi-Hopf algebra, and let $\yd ^{\rm fd}$ be
the category of finite dimensional left Yetter-Drinfeld modules
over $H$. If $M\in \yd ^{\rm fd}$, then $M^*\in \yd ^{\rm fd}$
(cf. \cite{bcp}). The action and coaction are given by
\begin{eqnarray}
&&\hspace*{-2cm}
(h\cd m^*)(m)=m^*(S(h)\cd m)\label{rdy1}\\
&&\hspace*{-2cm}\l _{M^*}(m^*)=\sum m^*_{(- 1)}\ot m^*_{(0)}=
\sum \limits _{i=1}^n \langle m^*, f^2\cd (g^1\cd {}_im)_{(0)}\rangle \nonumber\\
&&\hspace*{1cm}\smi (f^1(g^1\cd {}_im)_{(- 1)}g^2)\ot {}^im\label{rdy2}
\end{eqnarray}
for all $h\in H$, $m^*\in M^*$, $m\in M$. Here $f=\sum f^1\ot f^2$ is the twist
defined in (\ref{f}), $({}_im)_{i=\ov {1, n}}$ is a basis of $M$ and
$({}^im)_{i=\ov {1, n}}$ its dual basis. Moreover, $\yd ^{f. d.}$ is a rigid
monoidal category. For each object $M\in \yd ^{\rm fd}$,
the evaluation and coevaluation maps ($ev_M$ and $coev_M$, respectively)
are given by (\ref{qrig}).\\
In addition, if $B\in \yd ^{\rm fd}$ is a
Hopf algebra, then $B^*$ is a Hopf algebra in $\yd ^{f. d.}$. The structure
is the following.
\begin{itemize}
\item[-] The multiplication and unit are given by
\begin{eqnarray}
&&\hspace*{-2cm}
(\v * \psi )(b)=\langle \v , f^2\tilde{q}^2_2Y^3\smi
(\tilde{q}^1Y^1(p^1\cd b_{\un 2})_{(-1)}p^2)\cd b_{\un 1}\rangle \nonumber\\
&&\langle \psi , f^1\tilde{q}^2_1Y^2\cd (p^1\cd b_{\un 2})_{(0)}\rangle ,\label{dbm1}
\end{eqnarray}
\begin{equation}
1_{B^*}=\un{\va}\label{dbm2}
\end{equation}
for all $\v ,\psi \in B^*$, $b\in B$,
where $q_L=\sum \tilde{q}^1\ot \tilde{q}^2$ and $p_R=\sum p^1\ot p^2$ are the elements
defined in (\ref{ql}) and (\ref{qr}).
\item[-] the comultiplication and counit are given by the formulas
\begin{equation}
\un{\Delta }_{B^*}(\v )=\sum \limits _{i, j=1}^n
\langle \v , [(g^1\cd {}_jb)_{(-1)}g^2\cd {}_ib](g^1\cd {}_jb)_{(0)}\rangle 
{}^ib\ot {}^jb\label{dbc1}
\end{equation}
\begin{equation}
\un{\va }_{B^*}(\v )=\v (1_B),\label{dbc2}
\end{equation}
for any $\v \in B^*$, where $f^{-1}=\sum g^1\ot g^2$ was defined in (\ref{g}),
$({}_ib)_{i=\ov {1, n}}$ is a basis of $B$ and
$({}^ib)_{i=\ov {1, n}}$ the corresponding dual basis of $B^*$.
\item[-] the antipode is given by
\begin{equation}
\un{S}_{B^*}=\un{S}^*,
\mbox{${\;\;}$i. e. ${\;\;}$}
\un{S}_{B^*}(\v )=\v \circ \un{S},
\end{equation}
for all $\v \in B^*$.
\end{itemize}

\begin{proposition}\prlabel{5.4}
Let $B\in \yd ^{\rm fd}$ a Hopf algebra. Then $B^*$ is a right $B$-Hopf module, with 
structure:
\begin{equation}
\langle \v \leftharpoondown b, b^{'}\rangle =\sum \langle \v , [(U^1\cd b)_{(-1)}U^2\cd b^{'}]
\un{S}((U^1\cd b)_{(0)})\rangle ,\label{rhmbd1}
\end{equation}
\begin{equation}
\un{\r }_{B^*}(\v )=\sum \limits _{i=1}^n
(S(\tilde{p}^1)\cd {}_ib)_{(-1)}\cd [{}^ib*(\tilde{p}^2\cd \v )]
\ot (S(\tilde{p}^1)\cd {}_ib)_{(0)},\label{rhmbd2}
\end{equation}
for all $\v \in B^*$, $b, b^{'}\in B$, where
\begin{equation}\label{u}
U=\sum U^1\ot U^2:=\sum g^1S(q^2)\ot g^2S(q^1),
\end{equation}
$p_L=\sum \tilde{p}^1\ot \tilde{p}^2$,
$q_R=\sum q^1\ot q^2$ and $f^{-1}=\sum g^1\ot g^2$
are the elements defined by (\ref{ql}), (\ref{qr}) and (\ref{g}),
and $\{{}_ib\}_{i=\ov{1, n}}$ is a basis of $B$
with corresponding dual basis $\{{}^ib\}_{i=\ov{1, n}}$.
Moreover, 
$$B^{*{\rm co}B}=\{\Lambda \in B^*~|~ \sum (\tilde{p}^1\cd \v )*
(\tilde{p}^2\cd \Lambda) =\v (1_B)\Lambda \;\;for\;all\;\; \v \in B^*\}.$$
\end{proposition}

\begin{proof}
If $B$ is a Hopf algebra in a braided rigid monoidal category ${\mathcal C}$,
then $B^*$ is a right Hopf $B$-module, as follows.
\begin{itemize}
\item[-] the right $B$-module structure $\leftharpoondown : B^*\ot B\ra B$
on $B^*$ is the composition
\begin{equation}\label{rbmd}
\begin{array}{ccc}
B^*\ot B& \rTo^{l_{B^*\ot B}}&(B^*\ot B)\ot \un{1}\\
&\rTo^{(id_{B^*}\ot \un{S})\ot coev_B}&(B^*\ot B)\ot (B\ot B^*)\\
&\rTo^{a^{-1}_{B^*\ot B, B, B^*}}&((B^*\ot B)\ot B)\ot B^*\\
&\rTo^{a_{B^*, B, B}\ot id_{B^*}}&(B^*\ot (B\ot B))\ot B^*\\
&\rTo{(id_{B^*}\ot \un{m})\ot id_{B^*}}&(B^*\ot B)\ot B^*\\
&\rTo^{ev_B\ot id_{B^*}}&\un{1}\ot B^*\\
&\rTo^{r^{-1}_{B^*}}&B^*
\end{array}
\end{equation}
\item[-] the right $B$-comodule structure 
$\un{\r }_{B^*}:\ B^*\to B^*\ot B$
on $B^*$ is the composition
\begin{equation}\label{rbcmd}
\begin{array}{ccccc}
B^*&\rTo^{r_{B^*}}& \un{1}\ot B^*
&\rTo^{coev_B\ot id_{B^*}}&(B\ot B^*)\ot B^*\\
&\rTo^{a_{B, B^*, B^*}}&B\ot (B^*\ot B^*)
&\rTo^{id_B\ot \un{m}_{B^*}}&B\ot B^*\\
&\rTo^{c_{B, B^*}}&B^*\ot B.&&
\end{array}
\end{equation}
\end{itemize}
Let $\gamma =\sum \gamma ^1\ot \gamma ^2$ and $f^{-1}=\sum g^1\ot g^2$
be the elements defined in (\ref{gd}) and (\ref{g}). By (\ref{rbmd}),
we have,
for all $\v \in B^*$ and $b, b^{'}\in B$:
\begin{eqnarray*}
&&\hspace*{-1cm}
\langle \v \leftharpoondown b, b^{'}\rangle\\
 &=&\sum
\langle \v , S(X^1p^1_1)\a \cd [((X^2p^1_2\cd
\un{S}(b))_{(-1)}X^3p^2\cd b^{'})
(X^2p^1_2\cd \un{S}(b))_{(0)}]\rangle\\
{\rm (\ref{mal},\ref{y3})}&=&\sum
\langle \v , [(S(X^1p^1_1)_1\a _1X^2p^1_2\cd
\un{S}(b))_{(-1)}S(X^1p^1_1)_2\a _2X^3p^2\cd b^{'}]\\
&&\hspace*{1cm}(S(X^1p^1_1)_1\a _1X^2p^1_2\cd \un{S}(b))_{(0)}\rangle
 \end{eqnarray*}
\begin{eqnarray*}
{\rm (\ref{gdf},\ref{ca})}&=&\sum
\langle \v , [(g^1S(X^1_2p^1_{(1, 2)})\gamma ^1X^2p^1_2\cd
\un{S}(b))_{(-1)}g^2S(X^1_1p^1_{(1, 1)})\gamma ^2
X^3p^2\cd b^{'})]\\
&&\hspace*{1cm}(g^1S(X^1_2p^1_{(1, 2)})\gamma ^1X^2p^1_2\cd \un{S}(b))_{(0)}\rangle \\
{\rm (\ref{gd},\ref{q3},\ref{q5})}&=&\sum
\langle \v , g^1S(Y^2p^1_{(1, 2)})\a Y^3p^1_2\cd
\un{S}(b))_{(-1)}g^2S(Y^1p^1_{(1, 1)})\a p^2\cd b^{'}]\\
&&\hspace*{1cm}(g^1S(Y^2p^1_{(1, 2)})\a Y^3p^1_2\cd \un{S}(b))_{(0)}\rangle  \\
{\rm (\ref{q1},\ref{q5},\ref{q6})}&=&\sum
\langle \v , [(g^1S(Y^2)\a Y^3\cd \un{S}(b))_{(-1)}g^2S(Y^1)\cd b^{'}]\\
&&\hspace*{1cm}(g^1S(Y^2)\a Y^3\cd \un{S}(b))_{(0)}\rangle \\
{\rm (\ref{qr},\ref{u},\ref{smorf})}&=&\sum
\langle \v , [\un{S}(U^1\cd b)_{(-1)}U^2\cd b^{'}]\un{S}(U^1\cd b)_{(0)}\\
{\rm (\ref{smorf})}&=&\sum
\langle \v , [(U^1\cd b)_{(-1)}U^2\cd b^{'}]\un{S}((U^1\cd b)_{(0)})
\end{eqnarray*}
which is just (\ref{rhmbd1}). (\ref{rhmbd2}) follows easily by (\ref{rbcmd}),
the details are left to the reader. Finally, by (\ref{rbcmd})
we have
\begin{eqnarray*}
&&\hspace*{-2cm}\Lambda \in B^{*co(B)}
\Longleftrightarrow \un{\r}_{B^*}(\Lambda )=\Lambda \ot 1_B\\
&\Longleftrightarrow &c^{-1}_{B, B^*}\circ \un{\r }_{B^*}(\Lambda )=
c^{-1}_{B, B^*}(\Lambda \ot 1_B)\\
&\Longleftrightarrow &\sum \limits _{i=1}^n
S(\tilde{p}^1)\cd {}_ib\ot {}^ib*(\tilde{p}^2\cd \Lambda )=1_B\ot \Lambda\\
&\Longleftrightarrow &
\sum (\tilde{p}^1\cd \v )*(\tilde{p}^2\cd \Lambda )=\v (1_B)\Lambda,
~{\rm~for~all~}\v \in B^*.
\end{eqnarray*}
\end{proof}

We define the space of left integrals by $I_l(B^*)=B^{*co(B)}$. From the Fundamental
Theorem for Hopf modules, we then obtain.

\begin{corollary}\colabel{5.5}
Let $H$ be a quasi-Hopf algebra and $B$ a finite dimensional Hopf algebra
in $\yd $. Then $I_l(B^*)\ot B\simeq B^*$ as right $B$-Hopf modules. In particular,
$dim_k(I_l(B^*))=1$.
\end{corollary}

Now, let $H$ be a quasi-Hopf algebra and $H_0$ the $H$-module algebra
described in \seref{4}. If $(H,R)$ is quasitriangular, then $H_0$ is
a Hopf algebra in $\yd$, see \cite{bn}. The additional structure is the
following.
\begin{eqnarray}
&&
{\bf \lambda }_{H_0}(h)=\sum R^2\ot R^1\tr h,\label{scshz}\\
&&
\un{\Delta }(h)=\sum h_{\un{1}}\ot h_{\un{2}}\\
&&\hspace*{10mm}=\sum x^1X^1h_1g^1S(x^2R^2y^3X^3_2)
\ot x^3R^1\tr
y^1X^2h_2g^2S(y^2X^3_1),\label{und}\\
&&\un{\va }(h)=\va (h),\label{unva}\\
&&\un{S}(h)=\sum
X^1R^2p^2S(q^1(X^2R^1p^1\tr h)S(q^2)X^3),\label{unant}
\end{eqnarray}
for all $h\in H$, where $R=\sum R^1\ot R^2$ and
$f^{-1}=\sum g^1\ot g^2$, $p_R=\sum p^1\ot p^2$ and $q_R=\sum q^1\ot q^2$
are the elements defined by (\ref{g}), (\ref{qr}) and (\ref{qra}). 
By the above arguments, if $H$ is a finite
dimensional Hopf algebra, then
$H_0^*$ is also a Hopf algebra in $\yd $, with structure
\begin{eqnarray}
&&\hspace*{-2cm}(\v * \Psi )(h)=\sum \langle \v , f^2\ov{R}^2\tr h_{\un{1}}\rangle 
\langle \Psi , f^1\ov{R}^1\tr h_{\un{2}}\rangle \label{dmhz1}\\
&=&\sum \langle \v , f^2\tr Y^2\ov{R}^2X^1x^1_1h_1g^1S(Y^3x^3)\rangle \nonumber \\
&&\hspace*{1cm}
\langle \Psi , f^1Y^1\ov{R}^1\tr X^2x^1_2h_2g^2S(X^3x^2)\rangle ,\label{dmhz2}\\
&&\hspace*{-2cm}
1_{H_0^*}=\un{\va },\label{duhz}\\
&&\hspace*{-2cm}
\un{\Delta }_{H_0^*}(\v )=\sum \limits _{i, j=1}^n\langle \v , (R^2g^2\tr {}_ie)(R^1g^1\tr {}_je)\rangle 
{}^ie\ot {}^je\label{dcmhz},\\
&&\hspace*{-2cm}
\un{\va }_{H_0^*}(\v )=\v (\b ),\label{dchz}\\
&&\hspace*{-2cm}
\un{S}_{H_0^*}(\v) =v \circ \un{S},\label{danthz}
\end{eqnarray}
for all $h\in H$ and $\v \in H^*$, where $R^{-1}=\sum \ov{R}^1\ot \ov{R}^2$,
$\{{}_ie\}_{i=\ov{1, n}}$ is a basis of $H$ and
$\{{}^ie\}_{i=\ov{1, n}}$ the corresponding dual basis of $H^*$.
In this particular case we have
$$
I_l(H_0^*)=\{\Lambda \in H^*\mid
\sum \Lambda (S(\tilde{p}^2)f^1\ov{R}^1\tr h_{\un{2}})
S(\tilde{p}^1)f^2\ov{R}^2\tr h_{\un{1}}
=\Lambda (b)\b ,\;\;{\rm for\;all}\; h\in H\}.
$$


\end{document}